%% file: RealPolynomials2018.tex
\documentclass[12pt]{article}
\usepackage{amsmath}
\usepackage{amssymb}
\usepackage{enumitem}
\usepackage{graphicx}
\usepackage{color}
\usepackage{hyperref}

 \font\ncsc=cmcsc10
 \font\ntt=cmtt12

\def\CP{{{\mathbb C}{\rm P}}}

\def\C{\mathbb{C}}

\def\R{\mathbb{R}}

\def\qed{{\hfill $\Diamond$} \bigskip}
\def\Aut{{\rm Aut}}
\def\eps{\varepsilon}
\def\T{\mathbf{T}}
\def\W{\mathbf{W}}
\def\B{\mathbf{B}}

\def\Fo{F^{\rm odd}}
\def\Fe{F^{\rm even}}
\def\so{s^{\rm odd}}
\def\se{s^{\rm even}}
\def\base{S^2} 
\def\conj{\operatorname{conj}} 
\def\sou{\operatorname{source}}
\def\tar{\operatorname{target}}

\def\Gh{{\widehat{\Gamma}}}

\usepackage{theorem}

\newtheorem{theorem}{Theorem}

\newtheorem{proposition}{Proposition}[section]
\newtheorem{corollary}[theorem]{Corollary}
\newtheorem{lemma}[proposition]{Lemma}

{\theorembodyfont{\rmfamily}
\newtheorem{definition}[proposition]{Definition}
\newtheorem{example}[proposition]{Example}
\newtheorem{remark}[proposition]{Remark}
\newtheorem{notation}[proposition]{Notation}
}

\author{Ilia Itenberg, Dimitri Zvonkine}

\title{Hurwitz numbers for real polynomials}

\begin{document}

\maketitle

\begin{abstract}
We consider the problem of defining and computing real analogs of polynomial
Hurwitz numbers, in other words, the problem of counting properly normalized
real polynomials with fixed ramification profiles over real branch points. We
show that, provided the polynomials are counted with an appropriate sign, their
number does not depend on the order of the branch points on the real line. We
study generating series for the invariants thus obtained, determine necessary
and sufficient conditions for the vanishing and nonvanishing of these
generating series, and obtain a logarithmic asymptotic for the invariants as
the degree of the polynomials tends to infinity.
\end{abstract} 

\section{Introduction}

\subsection{Counting polynomials}
Let $P
\in \C[z]$ be a degree~$n$ polynomial with complex coefficients. 
To any $w \in \C$ we can assign a partition $\Lambda_w$ of~$n$ given by 
the orders of the roots of $P(z) = w$. 
This partition is called the {\em ramification type} of~$w$. The point $w$ is a branch point 
if and only if 
$\Lambda_w \not= (1, \dots ,1)$. If $w$ is a branch point we can consider the set of critical points $z$ such that $P(z)=w$. The multiplicities of these critical points, that is, the orders of vanishing of $P'(z)$ compose a partition $\lambda_w$. This partition is obtained from $\Lambda_w$ by subtracting~$1$ 
from every element of 
$\Lambda_w$
and eliminating the zeros. We 
call 
$\lambda_w$ 
the {\em reduced ramification type} of~$w$. 

The multiplicity of 
a branch point $w$ equals $|\lambda_w| = n - l(\Lambda_w)$,
where $|\lambda|$ is the sum of elements in $\lambda$, while $l(\Lambda)$ 
is the number of elements is~$\Lambda$. In particular, we have
$$
\sum_{w \in \C} |\lambda_w| = \sum_{w \in \C} (n - l(\Lambda_w)) = n-1.
$$
A branch point~$w$ is {\em simple} if it is of multiplicity~1, that is, $\Lambda_w= (2, 1, \dots, 1)$
and $\lambda_w = (1)$. 

We say that $P$ is {\em normalized} if it has the form   
$$
P(z) = z^n + a_2 z^{n-2} + \cdots + a_n.
$$ 

Counting the normalized polynomials with given branch points and their ramification types is a classical problem of enumerative geometry. It is equivalent to enumerating minimal factorizations of an $n$-cycle into a product of permutations of given cycle types in the symmetric group $S_n$. It is also equivalent to computing the so-called {\em polynomial Hurwitz numbers}. These numbers enumerate ramified coverings of the sphere by the sphere with one point of total ramification (corresponding to $\infty \in \CP^1$) and several other branch points with prescribed ramification profiles. The problem of counting normalized polynomials was posed 
by 
V.~Arnold (see~\cite{Arnold} Problem 1996-8), who also envisaged the possibility of studying the real case (Problem 1991-2). The problem for complex polynomials was completely solved in~\cite{GouJac, Zvonkine, LanZvo}. It is important to note that the answer does not depend on the positions of the branch points, but only on their ramification types.

Now assume that all the branch points are real. In 
this case it makes sense to count the {\em real} normalized polynomials with a given set of branch points and their ramification types. 
However, in this enumerative problem the answer 
in general does depend on the order of the branch points on the real line.
One of the goals 
of this paper is to show that the answer can be made invariant if we count each real polynomial with an appropriate sign.
Such a phenomenon was observed in various real enumerative problems; 
the first significant example is the Welschinger theorem \cite{Welschinger} providing an invariant signed count
of real rational curves in $4$-dimensional real symplectic manifolds. 

There exist other classical counting problems for ramified coverings that admit real analogs. For instance, the simplest problem is to count degree~$d$ ramified coverings of the sphere by the sphere with $2d-2$ simple branch points. The answer to this problem in the complex setting is $(2d-2)! \, d^{d-3}$; it is a particular case of a more general formula found by Hurwitz (see~\cite{Strehl}) and was rediscovered many times since then.
A real version of this counting problem in the case when all branch points are real was solved by B.~Shapiro and A.~Vainshtein in~\cite{ShaVai}.

The most general covering counting problem is to enumerate all ramified coverings of a genus $h$ surface by a genus $g$ surface with fixed ramification profiles over fixed branch points. Answers to this problem are usually called {\em Hurwitz numbers}. There is no closed formula for these numbers in either complex or real case, but some more or less practical methods to compute them. One approach uses the representation theory 
of the symmetric groups. We refer to Zagier's appendix to~\cite{LanZvoBook} for a review of the complex case and to A. Cadoret's work~\cite{Cadoret} for the real case (with $h = 0$). Another approach 
is based on a tropical correspondence theorem proved by B.~Bertrand, E.~Brugall\'e and G.~Mikhalkin. The complex case is treated in their paper~\cite{BBM}; for the real case, see the work by H.~Markwig and J.~Rau~\cite{MarRau}.

In the complex case, as well as in the setting of~\cite{ShaVai}, the invariance of Hurwitz numbers (that is, independence of the number of coverings from the positions of the branch points, provided that the ramification profiles are fixed) is immediate.
The papers \cite{Cadoret} and \cite{BBM} do not contain invariance statements: they enumerate all coverings 
sign `$+$' which, in general, gives rise to different Hurwitz numbers for different positions of branch points.
So far, our attempts to generalize the signed count and the invariance theorem to this general situation have failed. 

\subsection{The $s$-numbers}  
\begin{figure}[h]
\begin{center}
\ 
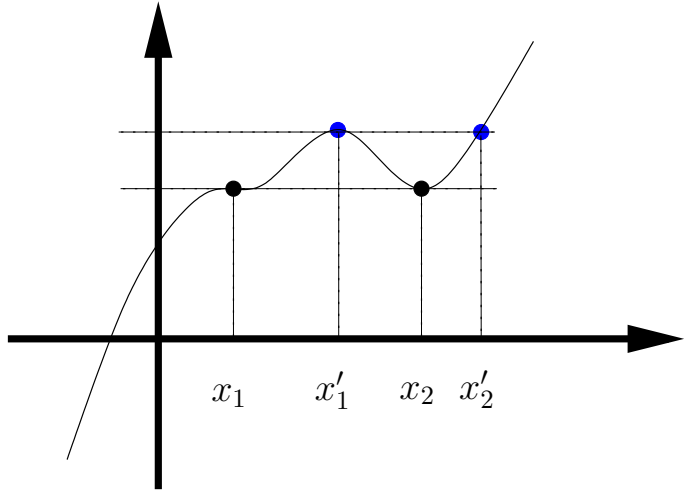
\caption{This polynomial has exactly two disorders: $(x_1, x_2)$ and $(x_1',x'_2)$} \label{Fig:disorder} 
\end{center}
\end{figure}

\begin{definition} \label{Def:disorderP}
Let $P \in \R[x]$ be a normalized real polynomial. 
A {\em disorder} of $P$ is a pair of real numbers $x_1 < x_2$ 
such that $P(x_1) = P(x_2)$ and the ramification order of $P$ at $x_1$ is greater than at~$x_2$, see Figure~\ref{Fig:disorder}.
\end{definition}

\begin{definition} \label{Def:signP}
Let $P \in \R[x]$ be a normalized real polynomial. The {\em sign} $\eps(P)$ of~$P$ is equal to $(-1)^d$, where $d$ is the total number of disorders in~$P$.

Given a real number $w$ we also define the {\em $w$-sign} of $P$ as $(-1)^{d_w}$, where $d_w$ is the number of disorders $x_1 < x_2$ of $P$ such that $P(x_1) = P(x_2) = w$. 
\end{definition} 

Let $k$ and $n$ be two positive integers, $k < n$. 
Choose a sequence $w_1, \dots, w_k$ 
of pairwise distinct real numbers, 
and let $(\Lambda_1, \dots, \Lambda_k)$ be a 
sequence of partitions of~$n$ such that 
$$
\sum_{i=1}^k (n - l(\Lambda_i)) = n-1.
$$ 
\begin{theorem}[the invariance theorem] \label{Thm:invariance}
Consider the set $S_{\Lambda_1, \ldots, \Lambda_k}(w_1, \ldots, w_k)$ 
of real normalized polynomials $P$ with branch points $w_1, \dots, w_k \in \R$ 
such that for every~$i$ the ramification type of $w_i$ is $\Lambda_i$. 
Then, the sum of signs 
$$
\sum_{P \in S_{\Lambda_1, \ldots, \Lambda_k}(w_1, \ldots, w_k)}\eps(P) 
$$ 
does not depend on $w_1, \dots, w_k$ {\rm (}in particular, on their order on the real line{\rm )\/}, 
but only on the partitions $\Lambda_1, \dots, \Lambda_k$.
\end{theorem} 

\begin{definition}
We call the sum
$$
\sum_{P \in S_{\Lambda_1, \ldots, \Lambda_k}(w_1, \ldots, w_k)} \eps(P) 
$$ 
from the previous theorem the {\it $s$-number} of real polynomials with given ramification type. 
\end{definition} 


\begin{example}
Let $n=4$ and $k = 2$. 
Put 
$\Lambda_1 = (2, 2)$ and $\Lambda_2 = (2, 1, 1)$. 
If $w_1 < w_2$, there are two real polynomials with ramification type  
$(\Lambda_1, \Lambda_2)$. 
Their graphs are shown in 
Figure \ref{Fig:two_solutions}. 
The polynomial $P_1$ has three real critical points, and its sign equals $\eps(P_1) = -1$. 
The polynomial $P_2$ has one real and two complex conjugate critical points, and its sign equals $\eps(P_2) = 1$. 
\begin{figure}[h]
\begin{center}
\ 
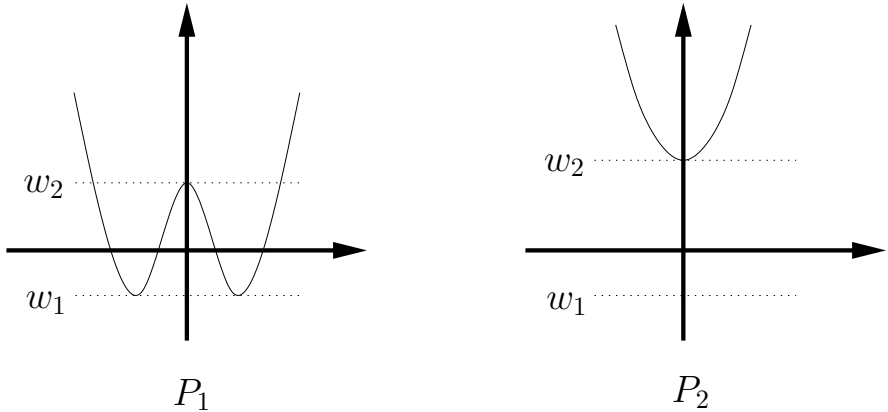
\caption{Two polynomials in the case $w_1 < w_2$} \label{Fig:two_solutions} 
\end{center}
\end{figure} 

If $w_1 > w_2$, there are no real polynomials with ramification type $(\Lambda_1, \Lambda_2)$. Thus, 
in both cases 
the $s$-number of real polynomials is equal to~$0$. 
\end{example} 

\subsection{Generating series}
In the second part of the paper we consider generating series
$F_{\lambda_1, \dots, \lambda_k}(q)$ for $s$-numbers. The coefficient of $q^m/m!$ in $F_{\lambda_1, \dots, \lambda_k}(q)$ is the $s$-number of polynomials with $k+m$ branch points. The first~$k$ points have reduced ramification types $\lambda_1, \dots, \lambda_k$, 
while the last~$m$ points are simple. 
Each generating series 
is decomposed into an even and an odd part: $\Fe$ 
enumerates real polynomials of even degrees, while $\Fo$ enumerates real polynomials of odd degrees. 

Let 
\begin{align*}
f(q) & = \tanh(q)= \frac{e^q - e^{-q}}{e^q+e^{-q}}= q - 2 \, \frac{q^3}{3!}+ 16 \, \frac{q^5}{5!} - 272 \, \frac{q^7}{7!} + \cdots,\\
g(q) & = \frac1{\cosh(q)}= \frac2{e^q+e^{-q}} = 1 - \frac{q^2}{2!}+ 5 \, \frac{q^4}{4!} - 61 \, \frac{q^6}{6!} + \cdots.
\end{align*}

\begin{theorem} \label{Thm:fg}
For any partitions $\lambda_1, \dots, \lambda_k$, the generating series
$$
\Fe_{\lambda_1, \dots, \lambda_k}(q)
$$
is a polynomial in $q$ and $f(q)$ with rational coefficients, while the generating series
$$
\Fo_{\lambda_1, \dots, \lambda_k}(q)
$$
is equal to $g(q)$ multiplied by a polynomial in $q$ and $f(q)$ with rational coefficients.
\end{theorem} 

\begin{figure}[h]
\begin{center}
\includegraphics[width=25em]{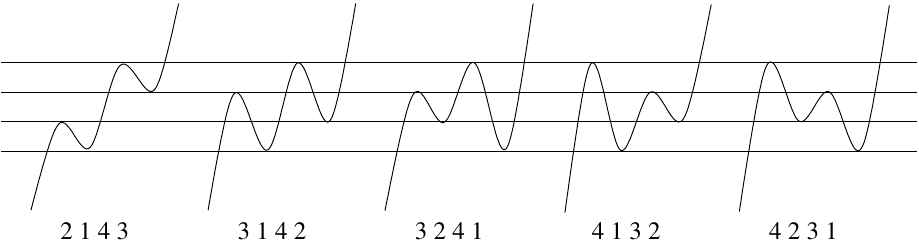}
\caption{The 5 real polynomials with 4 simple branch points and the corresponding alternating permutations. The number~5 is the coefficient of $q^4/4!$ in the power series~$g$.} \label{Fig:alternating}
\end{center}
\end{figure}

\begin{example}
We have $\Fe_\emptyset = f$ and $\Fo_\emptyset = g$. 
Indeed, when the set of partitions is empty, the polynomials only have simple critical values, so all critical points are real. 
If we number the critical points $x_1, \dots, x_{n-1}$ and the critical values $w_1, \dots, w_{n-1}$ in the increasing order on the real line, 
we obtain an alternating permutation~$\sigma$ given by $P(x_i) = w_{\sigma(i)}$, see Figure~\ref{Fig:alternating}.
The knowledge of $w_1, \dots, w_{n-1}$ and $\sigma$ determines the polynomial~$P$ uniquely (see Corollary~\ref{correspondence1} and Lemma~\ref{Lem:BaseToTree} for a generalization of this fact).
On the other hand, $f$ and $g$ are the generating series for the well-known Euler-Bernoulli numbers that enumerate alternating permutations, see the survey~\cite{Stanley} or the wikipedia entry~\url{Alternating_permutations}.
\end{example}

\begin{remark}
Note that the superscripts ``even'' and ``odd'' refer to the parity of~$n$, that is, the degree of~$P$. 
The series $\Fe$ and $\Fo$ may be even or odd depending on the partitions $\lambda_i$. 
More precisely, we have $m + \sum |\lambda_i| = n-1$, 
so that, for instance, the series $\Fe$ is even if $\sum |\lambda_i|$ is odd and odd if $\sum |\lambda_i|$ is even. 
\end{remark} 

\begin{theorem} \label{Thm:evenvanishing}
The generating series
$$
\Fe_{\lambda_1, \dots, \lambda_k}(q)
$$
is not identically~0
if and only if in each partition $\lambda_i$ every even number appears an even number of times 
and at most one odd number appears an odd number of times.
\end{theorem}

\begin{theorem} \label{Thm:oddvanishing}
The generating series
$$
\Fo_{\lambda_1, \dots, \lambda_k}(q)
$$
is not identically~0
if and only if in each partition $\lambda_i$ at most one odd number appears 
an odd number of times and at most one even number appears an odd number of times.
\end{theorem} 

The following statement concerns 
logarithmic asymptotics of $s$-numbers. 
Denote by $\se_{\lambda_1, \dots, \lambda_k}(m)$ and $\so_{\lambda_1, \dots, \lambda_k}(m)$ the coefficients of $q^m/m!$ in $\Fe_{\lambda_1, \dots, \lambda_k}(q)$ and $\Fo_{\lambda_1, \dots, \lambda_k}(q)$ respectively. These coefficients are, of course, the $s$-numbers of polynomials with given reduced ramification profile and $m$ 
additional simple branch points. 

\begin{theorem}\label{Thm:logasymp} 
Assume that the nonvanishing conditions of Theorem~\ref{Thm:evenvanishing} are satisfied. 
Then, we have 
$$
\ln |\se_{\lambda_1, \dots, \lambda_k}(m)|
\substack{\; \\ \; \\ \; \\ \sim\\ m \to \infty \\ m \;\; {\rm even}} m \ln m
\qquad \mbox{for} \; \, \sum |\lambda_i| \;\, \mbox{odd},
$$ 
$$
\ln |\se_{\lambda_1, \dots, \lambda_k}(m)|
\substack{\; \\ \; \\ \; \\ \sim\\ m \to \infty \\ m \;\; {\rm odd}} m \ln m
\qquad \mbox{for} \; \, \sum |\lambda_i| \;\, \mbox{even}.
$$ 
Assume that the nonvanishing conditions of Theorem~\ref{Thm:oddvanishing} are satisfied. 
Then, we have 
$$
\ln |\so_{\lambda_1, \dots, \lambda_k}(m)|
\substack{\; \\ \; \\ \; \\ \sim\\ m \to \infty \\ m \;\; {\rm odd}} m \ln m
\qquad \mbox{for} \; \, \sum |\lambda_i| \;\, \mbox{odd},
$$ 
$$
\ln |\so_{\lambda_1, \dots, \lambda_k}(m)|
\substack{\; \\ \; \\ \; \\ \sim\\ m \to \infty \\ m \;\; {\rm even}} m \ln m
\qquad \mbox{for} \; \, \sum |\lambda_i| \;\, \mbox{even}.
$$ 
\end{theorem} 

\begin{remark}
Consider the number of complex polynomials with $k$ branch points of reduced ramification types $\lambda_1, \dots, \lambda_k$ and $m$ additional simple branch points. It is easy to deduce from the explicit formulas of~\cite{Zvonkine} that the logarithmic asymptotics for this number as $m \to \infty$ is $m \ln m$. Thus, when the nonvanishing conditions of Theorems~\ref{Thm:evenvanishing} and~\ref{Thm:oddvanishing} are satisfied, 
the number of complex polynomials, 
the absolute value of the $s$-number of real polynomials, and the actual number of real polynomials, which lies between the two latter numbers, all have the same logarithmic asymptotic. A similar phenomenon has been observed and proved in many other situations, see, for instance, \cite{IKS1, IKS2}.  
\end{remark}

\begin{remark}
In the complex case the difference between counting ramified coverings of the sphere and normalized polynomials is rather trivial: the answers differ by a factor of $n$ due to the the fact that the change of variables $z \mapsto \sqrt[n]{1}\, z$ changes the normalized polynomial, but not the ramified covering. In the real case, however, the difference is more subtle, since a real polynomial can be normalized by a real change of variables in 0, 1, or 2 ways depending on its parity and the sign of the leading coefficient. It seems that the problem that admits a nice real version is the counting of normalized polynomials rather than ramified coverings. 
\end{remark} 

\begin{remark}
Many of our results, in particular Theorems~\ref{Thm:invariance} and~\ref{Thm:fg}, remain true if we allow the real polynomials to have pairs of complex conjugate branch points with equal ramification types. The notion of a real polynomial dessin, that we introduce in Section~\ref{dessin} as the preimage under~$P$ of the real axis, must then by replaced by the preimage of a connected graph 
containing the real axis and the complex branch points. 
The description of such generalized dessins makes the proofs rather cumbersome without adding much to the understanding; therefore we chose to restrict ourselves to the case of real branch points.
\end{remark} 

\paragraph{Plan of the paper.} In Section~\ref{dessin} we introduce real polynomial dessins. Such dessins are a strandard way to capture the combinatorial structure of the preimage $P^{-1}(\R)$. They are in a one-to-one correspondence with real normalized polynomials. In Section~\ref{Sec:trees} we prove a combinatorial theorem on black and white real trees. This theorem is equivalent to  Theorem~\ref{Thm:invariance} (the invariance theorem) in the particular case of two critical values. In Section~\ref{proof} we deduce the full statement of the Theorem~\ref{Thm:invariance} from this particular case. Finally, in Section~\ref{generating} we study the generating series for $s$-numbers of real polynomials and prove Theorems~\ref{Thm:fg}, \ref{Thm:evenvanishing}, \ref{Thm:oddvanishing}, and~\ref{Thm:logasymp}.

\section{Real polynomial dessins}\label{dessin}
Real polynomial dessins were first introduced and used to enumerate real polynomials by S.~Barannikov in~\cite{Barannikov}.

From now and till the end of this section we fix two positive integers $k < n$ and a sequence $(\Lambda_1, \dots, \Lambda_k)$ of partitions of~$n$ such that 
$$
\sum_{i=1}^k (n - l(\Lambda_i)) = n-1.
$$

Let $c: S^2 \to S^2$ be an orientation reversing involution of a $2$-dimensional sphere $S^2$. We assume that the fixed point set $E \subset S^2$ of $c$ is homeomorphic to a circle and choose an orientation of this circle. Let us introduce a graph that captures the combinatorial structure of $P^{-1}(\R)$ for a polynomial~$P$.

\begin{definition} \label{Def:dessin}
A \emph{real polynomial dessin} of degree $n$ and type $(\Lambda_1, \ldots, \Lambda_k)$ in $\base$  
is an oriented graph $\Gamma\subset\base$ whose vertices are labelled
by elements of the set $\{ 1, 2, \dots, k, \infty \}$ in such a way that 
\begin{itemize}
\item the oriented graph $\Gamma$ (together with the labeling of the vertices)
is invariant under $c$; 
\item the circle $E$ is a union of edges of $\Gamma$;
\item
exactly one vertex of $\Gamma$ is labelled by $\infty$, and the degree of this vertex is $2n$; 
\item for each integer $1 \leq i \leq k$, the graph $\Gamma$ has exactly $l(\Lambda_i)$ vertices labelled by $i$ 
and their degrees are equal to the elements of $\Lambda_i$, multiplied by $2$; 
\item each edge of $\Gamma$ is one of the following $k + 1$ types:
$$\infty \to 1, \; 1 \to 2, \;  2 \to 3, \; \ldots, \; k - 1 \to k, \; k \to \infty,$$ 
where $i \to j$ means that the edge starts at a vertex labelled by $i$
and finishes at a vertex labelled by $j$; 
\item for any connected component $C$ of $\base \setminus \Gamma$,
each type of edges appears exactly once in the boundary $\partial C$ of $C$
(in particular, the orientation of the edges in $\partial C$
extends to an orientation of $C$).  
\end{itemize} 
\end{definition}

\begin{definition} \label{Def:DessinHomeo}
Two real polynomial dessins $\Gamma_1\subset S^2_1$ and $\Gamma_2 \subset S^2_2$
of degree $n$ and type $(\Lambda_1, \ldots, \Lambda_k)$
are {\em homeomorphic} if there exists a homeomorphism $\varphi: S^2_1 \to S^2_2$ such that
\begin{itemize}
\item $\varphi \circ c_1 = c_2 \circ \varphi$,
\item $\varphi$ respects the chosen orientations
of $E_1$ and $E_2$ (that is, $\varphi_{| E_1}: E_1 \to E_2$ is of degree $1$ with respect to the orientations
of $E_1$ and $E_2$),  
\item $\varphi(\Gamma_1) = \Gamma_2$, and $\varphi$ preserves the labels of vertices and the types of edges. 
\end{itemize} 
\end{definition}

The vertices of a real polynomial dessin $\Gamma \subset S^2$ that belong to $E$ are called {\it real}. 
The complement of the $\infty$-vertex in $E$ is a totally ordered set that can be identified with 
the real line $\R$. 
 
\begin{definition}\label{Def:DessinDisorder}
Consider a pair of real vertices $v_1$ and $v_2$ of $\Gamma$ that are labeled by the same number different from $\infty$.
We say that the pair formed by $v_1$ and $v_2$ is a {\it disorder} of $\Gamma$ if 
\begin{itemize}
\item 
$v_1$ is smaller than $v_2$
in $E \setminus \{\infty\}$,
\item the degree of $v_1$ is bigger than the degree of $v_2$. 
\end{itemize} 
\end{definition} 

\begin{definition}\label{Def:DessinSign}
The {\it sign} $\varepsilon(\Gamma)$ of a real polynomial dessin $\Gamma \subset S^2$
is $(-1)^{d(\Gamma)}$, where $d(\Gamma)$ is the number of disorders of $\Gamma$. 
\end{definition}

\begin{definition} \label{Def:increasing}
A real polynomial dessin $\Gamma \subset S^2$ is said to be {\it increasing}, if $E$ contains an edge of $\Gamma$
of type $k \to \infty$ 
such that 
the orientation of this edge coincides with the orientation of $E$.  
\end{definition} 

\begin{figure}
\begin{center}
\includegraphics[width=\textwidth]{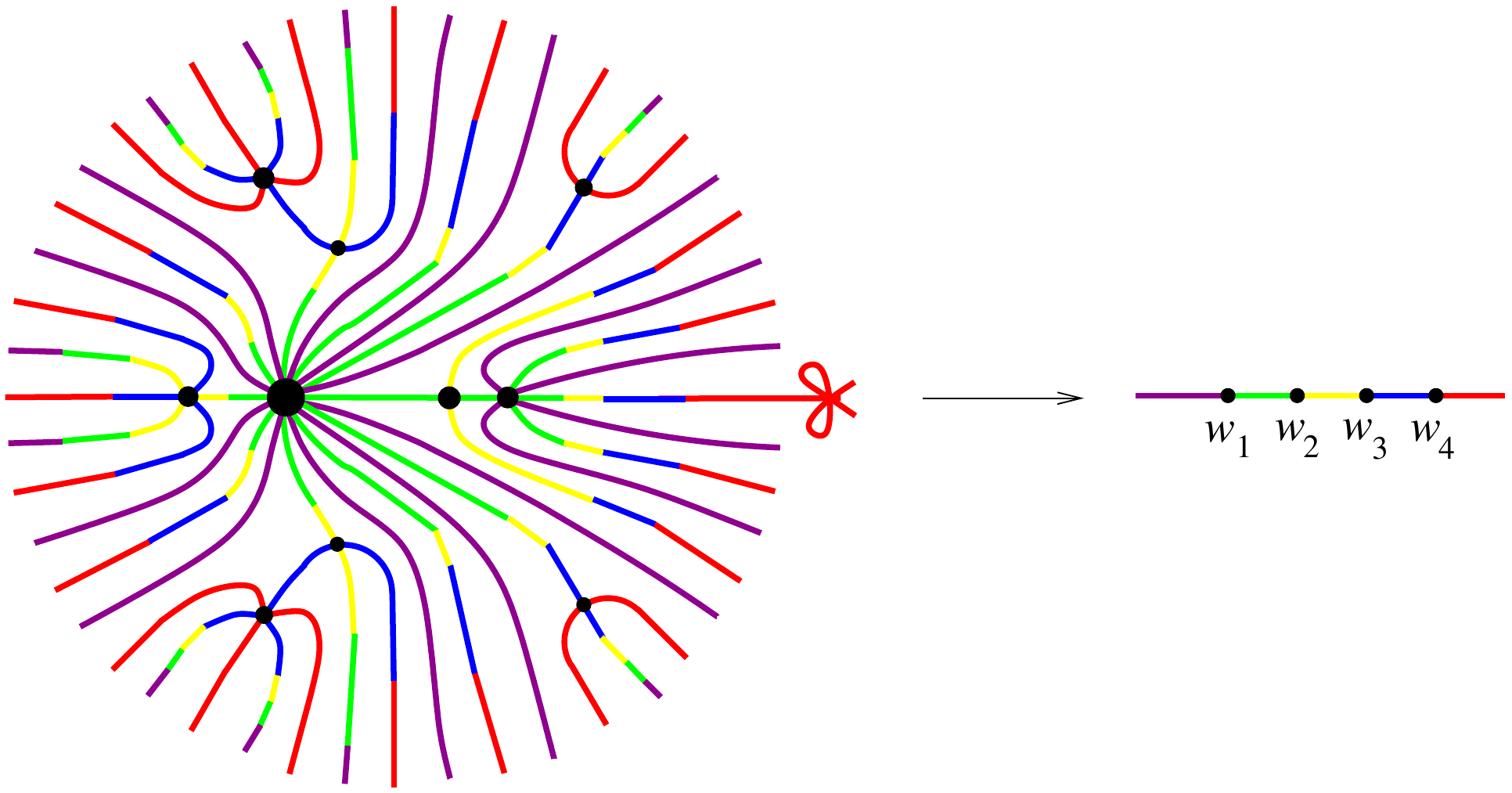}
\caption{An increasing real polynomial dessin obtained as the pull-back of the real axis under 
a real polynomial function.  
The edge of type $k \to \infty$ is marked with a ribbon.} \label{Fig:affinedessin}
\end{center}
\end{figure}

In the figures we 
represent the {\em affine dessins}, that is, the dessins without the unique vertex labeled by $\infty$,
see Figure~\ref{Fig:affinedessin}. 

A typical example of a real polynomial dessin is as follows.  
Let $\C P^1_{\sou} = \C_{\sou} \cup \infty_{\sou}$ and $\C P^1_{\tar} = \C_{\tar} \cup \infty_{\tar}$
be two copies of the complex projective line $\C P^1$. 
A holomorphic function 
$$P: \C P^1_{\sou} \to \C P^1_{\tar}$$ 
is {\it polynomial},  
if $P^{-1}(\infty_{\tar}) = \infty_{\sou}$. 
In affine coordinates in $\C_{\sou}$ and $\C_{\tar}$ such a function is indeed given by a polynomial.
A polynomial function $P: \C P^1_{\sou} \to \C P^1_{\tar}$ is {\it real} if 
$\conj_{\tar} \circ P = P \circ \conj_{\sou}$, where $\conj: \C P^1 \to \C P^1$ is the standard complex conjugation
in $\C P^1$. In the affine coordinate this means that $\overline{P}(z) = P(\overline{z})$, 
{\it i.e.}, 
that 
$P$ has real coefficients. 

Let $P: \C P^1_{\sou} \to \C P^1_{\tar}$ be a real polynomial function of degree $n$ and such that all the branch points $w_1, \ldots, w_k$ of $P$ in $\C_{\tar}$ are real. Assume that $w_1 < w_2 < \ldots < w_k$, 
and denote the ramification type of $P$ at $w_i$ by $\Lambda_i$, $i = 1$, $\ldots$, $k$. The polynomial function $P$ defines a real polynomial dessin $\Gamma_P \subset \C P^1_{\sou}$.  
As a set, $\Gamma_P$ is the pull-back $P^{-1}(\R P^1_{\tar})$, where $\R P^1_{\tar} \subset C P^1_{\tar}$ is the fixed point set of $\conj_{\tar}$.
The real polynomial dessin structure on~$\Gamma_P$ is introduced
as follows: 
the involution $c$ is $\conj_{\sou}$, the orientation of the fixed point set $\R P^1_{\sou}$ of $\conj_{\sou}$ is {\it positive} 
({\it i.e.}, induced by the order of $\R \subset \R P^1_{\sou}$), 
the only vertex of $\Gamma$ labelled by $\infty$ is $\infty_{\sou}$, the pull-backs of $w_1$, $\ldots$, $w_k$ 
are labelled by $1$, $\ldots$, $k$, respectively; 
the orientation of~$\Gamma_P$ is that induced from the positive
orientation of~$\R P^1_{\tar}$.  

\begin{proposition}\label{poly->dessin}
Let 
$P: \C P^1_{\sou} \to \C P^1_{\tar}$ 
be a real polynomial function of degree $n$ 
such that all the branch points $w_1, \ldots, w_k$ of $P$ in $\C_{\tar}$ are real.
Assume that $w_1 < w_2 < \ldots < w_k$, 
and denote the ramification type of $P$ at $w_i$ by $\Lambda_i$, $i = 1$, $\ldots$, $k$.
Then $\Gamma_P$ is a real polynomial dessin of degree $n$ and type $(\Lambda_1, \ldots, \Lambda_k)$.  
\end{proposition}

\paragraph{Proof.} Straightforward.  \qed

\vskip10pt

A real polynomial dessin corresponding to a real polynomial function $P$
is increasing if and only if $P(x) \to +\infty$ as $x \to +\infty$, 
in other words, if and only if the leading coefficient of $P$ is positive. 

A real polynomial function whose critical values are all real is said to be {\it totally real}. 
The following proposition shows that any real polynomial dessin can be identified with $\Gamma_P$
for a certain totally real polynomial function $P$.    

\begin{proposition}\label{dessin->poly}
Let 
$\Gamma \subset \base$ be a real polynomial dessin of degree $n$ 
and type $(\Lambda_1, \ldots, \Lambda_k)$, 
and let $w_1 < w_2 < \ldots < w_k$ be 
real numbers. 
Then, there exists  
a totally real polynomial function $P: \C P^1_{\sou} \to \C P^1_{\tar}$ 
of degree $n$ and type $(\Lambda_1, \ldots, \Lambda_k)$ 
such that $\Gamma_P$ is homeomorphic to $\Gamma$, 
and $w_i$ is a branch point of type $\Lambda_i$ of $P$ for each $i = 1$, $\ldots$, $k$. 
\end{proposition} 

\paragraph{Proof.} 
The points $w_1$, $\ldots$, $w_k$, and $\infty_{\tar}$ divide $\R P^1_{\tar}$ into $k + 1$ segments, 
which are called {\it non-critical}.
Denote by $S^+$ and $S^-$ the two semi-spheres of $\base$ which have the common boundary $E \subset \base$,
and construct a ramified $n$-fold covering $\Phi: \base \to \C P^1_{\tar}$ such that 
$\conj_{\tar} \circ \Phi = \Phi \circ c$
in the following way. 

Put the image under $\Phi$ of the point labelled by $\infty$ to be $\infty_{\tar}$. 
For each vertex $v \in S^+$ of $\Gamma$ such that $v$ is labelled by $i \in \{ 1, \dots, k \}$ set the image under $\Phi$ of $v$ to be $w_i$.
For each edge $e \subset S^+$ of $\Gamma$, 
let $\Phi$ send $e$ homeomorphically to the non-critical segment between the critical values corresponding to the extremal points of $e$  
(in such a way that the orientation of the edge corresponds to the positive orientation of $\R P^1_{\tar}$). Since the degree of every vertex of the dessin $\Gamma$ is even, the connected components of $\base \setminus \Gamma$ possess a chessboard coloring: for any two neighboring components, one of them is black, and the other is white.   
Extend $\Phi$ to $S^+$ sending each connected component of $\base \setminus \Gamma$
homeomorphically to one of the halves $\C P^1_{\tar} \setminus \R P ^1_{\tar}$ 
in such a way that all the connected components of the same color
are sent to the same half, and connected components of different colors
are sent to different halves.  
Finally, put $\Phi_{|S^-} = \conj_{\tar} \circ \Phi_{|S^+} \circ c$. 

The resulting ramified covering $\Phi: \base \to \C P^1_{\tar}$
is equivariant 
({\it i.e.}, satisfies the condition $\conj_{\tar} \circ \Phi = \Phi \circ c$) 
and allows one to lift the complex structure from $\C P^1_{\tar}$ to $\base$.
Namely, due to the Riemann existence theorem,
there exists a complex structure on $\base$ such that the map $\Phi: \base \to C P^1_{\tar}$
becomes holomorphic and the involution $c$ becomes anti-holomorphic.  
The uniqueness of the complex structure
on a $2$-dimensional sphere implies 
the existence of an equivariant biholomorphic isomorphism $\varphi: S^2 \to \C P^1_{\sou}$;
composing, if necessary, $\varphi$ with the multiplication by $-1$ in $\C P^1_{\sou}$,
we can assume that $\varphi$ respects the orientations of $E$ and $\R P^1_{\sou}$.
By construction, the map $\Phi \circ \varphi^{-1}$ is a totally real polynomial function of degree $n$ 
and ramification type ($\Lambda_1$, $\ldots$, $\Lambda_k$),
and the real polynomial dessins $\Gamma$ and $\Gamma_P$ are homeomorphic. 
\qed 

\vskip5pt 

A {\it real affine transformation} of $\C P^1_{\sou}$ is a homography of $\C P^1_{\sou}$ 
which preserves $\infty_{\sou}$ and commutes with $\conj_{\sou}$. In the affine coordinate it has the form $z \mapsto az+b$, with $a \in \R^*$ and $b \in \R$. Such a transformation is called {\it positive} if it respects the orientation of $\R P^1_{\sou}$ (in other words, if $a$ is positive) and {\em negative} otherwise.  Denote by $K^+_{\sou}$ the group 
of positive real affine transformations of $\C P^1_{\sou}$. 

\begin{corollary}\label{correspondence}
Let $w_1 < \ldots < w_k$ be real numbers. 
Then, the correspondence $P \mapsto \Gamma_P$
establishes a bijection between, on the one hand, the $K^+_{\sou}$-orbits on the set of totally real polynomial functions of degree $n$ 
with branch points $w_1$,  $\ldots$, $w_k$ of types $\Lambda_1$, $\ldots$, $\Lambda_k$, respectively, and, on the other hand,
the set of homeomorphism classes of real polynomial dessins of degree $n$ 
and type $(\Lambda_1, \ldots, \Lambda_k)$. 
\end{corollary} 

\paragraph{Proof.} The statement immediately follows from Proposition \ref{dessin->poly} and
the fact that the isomorphism $\varphi$ in the proof of Proposition \ref{dessin->poly} is unique up to
the action (by composition) of $K^+_{\sou}$
and $\conj_{\sou}$. 
\qed 

\begin{corollary}\label{correspondence1}
Let $w_1 < \ldots < w_k$ be real numbers. 
Then, the correspondence $P \mapsto \Gamma_P$ 
establishes a bijection between, on the one hand, the set $S_{\Lambda_1, \ldots, \Lambda_k}(w_1, \ldots, w_k)$  
of real normalized polynomials of degree $n$
with branch points $w_1$,  $\ldots$, $w_k$ of types $\Lambda_1$, $\ldots$, $\Lambda_k$, respectively, 
and, on the other hand, the set of homeomorphism classes 
of increasing {\rm (}see Definition~\ref{Def:increasing}{\rm )} real polynomial dessins of degree $n$ 
and type $(\Lambda_1, \ldots, \Lambda_k)$.
Furthermore, for any real normalized polynomial $P \in S$, the sign $\varepsilon(P)$ of $P$ coincides with $\varepsilon(\Gamma_P)$.  
\end{corollary} 

\paragraph{Proof.} A real polynomial $P$ can be normalized 
by a positive real affine transformation (an element of $K^+_{\sou}$) if and only if its leading coefficient 
is positive or, in other words, if and only if the corresponding dessin is increasing. In this case $P$ can be normalized in a unique way. 
Thus, Corollary \ref{correspondence} implies the first statement.
The statement about signs is immediate. 
\qed

\vskip10pt

Denote by $D^+_{\Lambda_1, \ldots, \Lambda_k}$ the set of homeomorphism classes
of increasing real polynomial dessins of degree $n$ and type $(\Lambda_1, \ldots, \Lambda_k)$.
We obtain the following statement. 

\begin{corollary}\label{signed_equality}
Let $w_1 < \ldots < w_k$ be real numbers.
Then, 
the $s$-number $\sum_{P \in S_{\Lambda_1, \ldots, \Lambda_k}(\omega_1, \ldots, \omega_k)}\varepsilon(P)$
is equal to $\sum_{\Gamma \in D^+_{\Lambda_1, \ldots, \Lambda_k}}\varepsilon(\Gamma)$. 
\qed 
\end{corollary}

The number $\sum_{\Gamma \in D^+_{\Lambda_1, \ldots, \Lambda_k}}\varepsilon(\Gamma)$ is called
the {\it $s$-number of increasing real polynomial dessins of degree $n$ and type $(\Lambda_1, \ldots, \Lambda_k)$}.

\section{Black and white trees} \label{Sec:trees}

In this section, we introduce auxiliary combinatorial objects that are used in the proof of the invariance theorem. They are equivalent to real polynomial dessins with $k=2$ critical values. 

\begin{definition} \label{Def:tree}
A {\em black and white tree} is a tree embedded into $\C$
whose vertices are colored in black and white in alternation.
A black and white tree is said to be {\it real} if it 
is invariant (including the colors) under the complex conjugation. 
Two real black and white 
trees are isomorphic if one can be transformed into the other by
an equivariant (with respect to the complex conjugation) homeomorphism of~$\C$ 
preserving the orientation of the real axis. 
\end{definition} 

\begin{remark}\label{remark-completion} 
To transform a real black and white tree into a real polynomial dessin
we add one vertex labeled~$\infty$ and, at each other vertex, insert an edge leading to the $\infty$-vertex between 
each pair of successive edges in such a way that the result is invariant under the complex conjugation. 
Thus, the degree of each vertex of the tree is multiplied by~2, while the degree of the $\infty$-vertex
is twice the number of edges of the tree. 
\end{remark} 

For shortness, most of the time in this section we  
just say {\em tree} instead of ``real black and white 
tree''. We always consider trees up to isomorphism. 


\begin{definition} \label{Def:realpart}
To a tree~$T$ we assign its {\em real part sequence}: the sequence of colors and degrees of its vertices lying on the real axis from left to right. 
The vertices of the real part $\R T = T \cap \R$ of $T$ are called 
{\em real vertices} of~$T$. The first and last real vertices are called the {\em border vertices}; 
the other real vertices are called {\em interior real vertices}.
\end{definition} 

\begin{definition}
A {\em disorder} of a tree is a pair of its real vertices of the same color such that the degree of the first one is greater
than the degree of the second one.
\end{definition} 

\begin{definition} \label{Def:sign}
The {\em sign} of a tree $T$ is 
$$
\eps(T) = (-1)^d,
$$
where $d$ is the total number of disorders in~$T$.
\end{definition}

\begin{definition} \label{Def:side}
A tree is {\em a white side tree} (respectively {\em a black side tree}) if its rightmost vertex on the real axis is white (respectively, black).
\end{definition}

\begin{example}
Up to isomorphism, 
there are exactly 12 black and white real trees with 4 edges. 
They are shown in 
Figure \ref{Fig:4edges} 
with their signs. The black side trees are on the left, while the white side trees are on the right. 

\begin{figure}
\begin{center}
\includegraphics[width=20em]{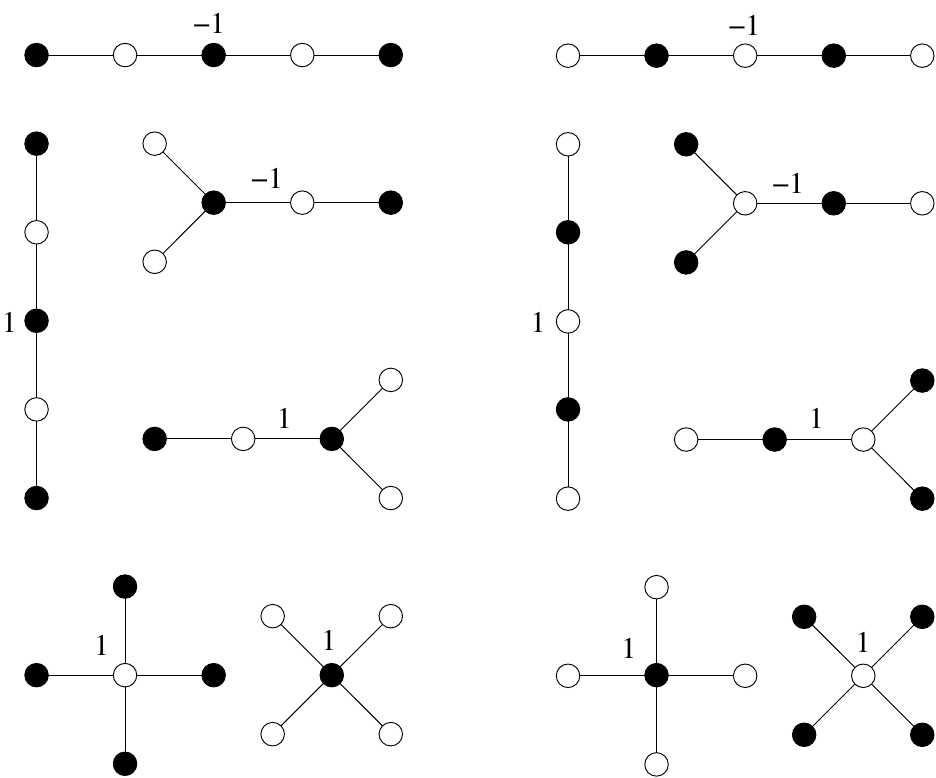} 
\caption{Black and white trees with $4$ edges} \label{Fig:4edges} 
\end{center}
\end{figure} 
\end{example} 

\begin{notation}
Fix a positive integer $n$ and two partitions $\Lambda_w$ and $\Lambda_b$ of~$n$. We will denote by $\T_{\Lambda_w,\Lambda_b}$, $\W_{\Lambda_w,\Lambda_b}$, and $\B_{\Lambda_w,\Lambda_b}$ the sets of trees, 
white side trees and black side trees, respectively,
whose degrees of white and black vertices are prescribed by $\Lambda_w$ and $\Lambda_b$. 
Thus, 
$$
\T_{\Lambda_w,\Lambda_b} = \W_{\Lambda_w,\Lambda_b} \sqcup \B_{\Lambda_w,\Lambda_b}.
$$
\end{notation}

\begin{theorem}\label{tree-invariance} 
Fix a positive integer $n$ and two partitions $\Lambda_w$ and $\Lambda_b$ of~$n$. We have
$$
\sum_{T \in \W_{\Lambda_w,\Lambda_b}} \!\!\!\!\!\!\! \eps(T) 
\; - \!\!\!\!
\sum_{T \in \B_{\Lambda_w,\Lambda_b}} \!\!\!\!\!\! \eps(T) =0.
$$
In other words, there is the same number of white side and black side trees 
whose degrees of white and black vertices are prescribed by $\Lambda_w$ and $\Lambda_b$
provided that we count the trees with signs given in Definition~\ref{Def:sign}.
\end{theorem} 

\begin{figure}[h]
\begin{center}
\includegraphics[width=25em]{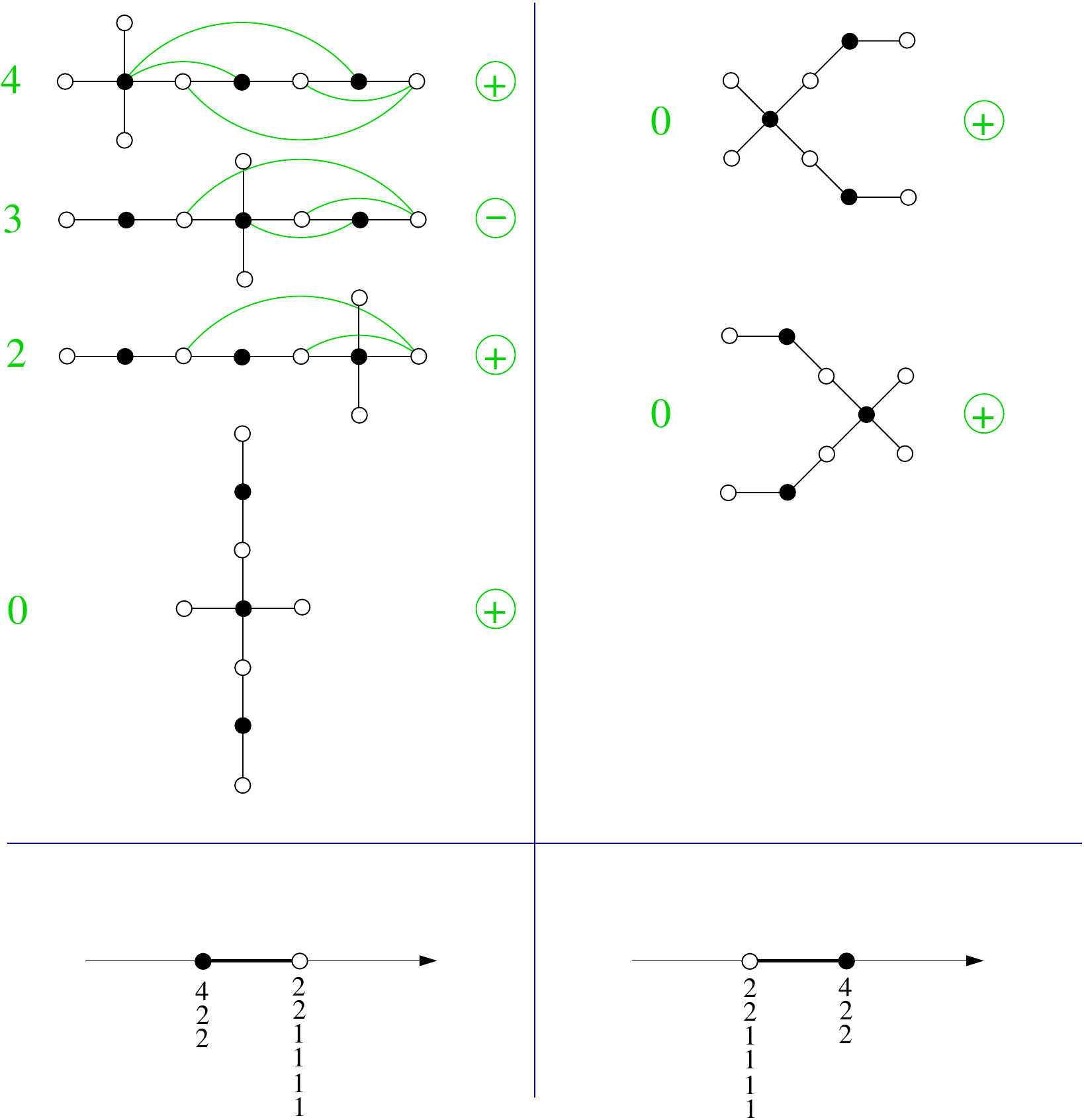}
\caption{White side trees are on the left, black side trees on the right. All trees have $\Lambda_b = (4,2,2)$, $\Lambda_w = (2,2,1,1,1,1)$. The disorders, their number, and the resulting signs are shown in green.} \label{Fig:tree-invariance}
\end{center}
\end{figure}

\begin{example}
Figure~\ref{Fig:tree-invariance} shows all trees with $\Lambda_b = (4,2,2)$, $\Lambda_w = (2,2,1,1,1,1)$. The disorders are shown as green arcs, their number is written to the left of each tree and the resulting sign to the right of each tree. The sum of signs is equal to 2 both for the white side trees on the left and for the black side trees on the right.
\end{example}

To prove the theorem we 
need an auxiliary way of weighting trees. 

\begin{definition} \label{Def:weight}
We define the {\em weight} $\omega(T)$ of a tree $T$ as follows.
If the real part sequence of a tree~$T$ is not symmetric, we have $\omega(T) = 0$. If $T$ has only one real vertex, then $\omega(T) = 1$. 
If the real part sequence of~$T$ is symmetric, $T$ has more than one real vertex, and the middle real vertex has the same color as the border vertices, then $\omega(T) = -1$. 
If the real part sequence of~$T$ is symmetric, $T$ has more than one real vertex, and the color of the middle real vertex is not the same as the color of the border vertices, then $\omega(T) = 1$.
\end{definition}

The statement of Theorem~\ref{tree-invariance} 
is a consequence of the following lemmas. 

\begin{lemma} \label{Lem:signweight}
Fix a positive integer $n$ and two partitions $\Lambda_w$ and $\Lambda_b$ of~$n$. Then we have 
$$
\sum_{T \in \W_{\Lambda_w,\Lambda_b}} \!\!\!\!\!\!\! \eps(T) \; - \!\!\!\!
\sum_{T \in \B_{\Lambda_w,\Lambda_b}} \!\!\!\!\!\! \eps(T) \; \; = \!\!
\sum_{T \in \W_{\Lambda_w,\Lambda_b}} \!\!\!\!\!\!\! \omega(T) \; - \!\!\!\! \sum_{T \in \B_{\Lambda_w,\Lambda_b}} \!\!\!\!\!\! \omega(T).
$$
In other words, the difference between the numbers of white and black side trees is the same whether we count the trees with signs of Definition~\ref{Def:sign} or with weights of Definition~\ref{Def:weight}.
\end{lemma}

\begin{figure}[h]
\begin{center}
\includegraphics[width=25em]{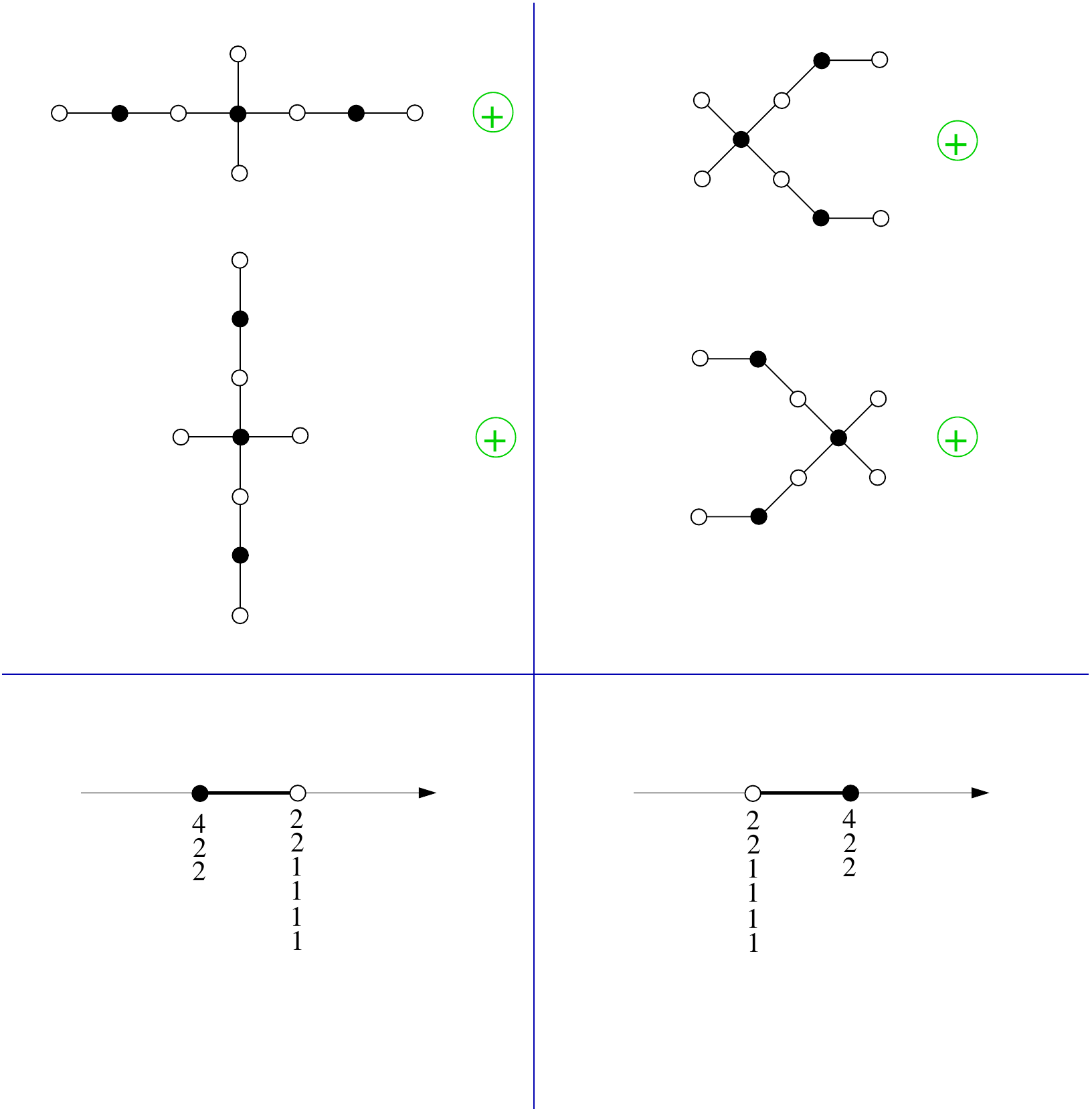}
\caption{White side trees are on the left, black side trees on the right. All trees have $\Lambda_b = (4,2,2)$, $\Lambda_w = (2,2,1,1,1,1)$. The weights of the trees are shown in green.} \label{Fig:signtoweight}
\end{center}
\end{figure}

\begin{example}
Figure~\ref{Fig:signtoweight} shows all trees with $\Lambda_b = (4,2,2)$, $\Lambda_w = (2,2,1,1,1,1)$ 
and symmetric real part sequence, together with their weights. 
Note that the three first trees on the white side of Figure~\ref{Fig:tree-invariance} 
got replaced by just one tree whose {\em weight} is equal to the sum of {\em signs} of the three trees.
\end{example} 

\begin{lemma} \label{Lem:weightzero}
We have
$$
\sum_{T \in \W_{\Lambda_w,\Lambda_b}} \!\!\!\!\!\!\! \omega(T) \; - \!\!\!\! \sum_{T \in \B_{\Lambda_w,\Lambda_b}} \!\!\!\!\!\! \omega(T) =0.
$$
\end{lemma}

\paragraph{Proof of Lemma~\ref{Lem:signweight}.} We start with the sum
$$
\sum_{T \in \W_{\Lambda_w,\Lambda_b}} \!\!\!\!\!\!\! \eps(T) \; - \!\!\!\!
\sum_{T \in \B_{\Lambda_w,\Lambda_b}} \!\!\!\!\!\! \eps(T)
$$
and reduce it by finding pairs of trees that cancel with each other.

Consider a tree $T \in \T_{\Lambda_w,\Lambda_b}$.

{\bf A. Symmetrizing the border vertices.} If the border vertices have the same color but different degrees, we can construct a new tree~$T'$ by interchanging them together with their forests. We have $\eps(T') = -\eps(T)$. This is due to the fact that the degrees of the border vertices are odd and therefore are never equal to the degrees of the other real vertices that are even. The trees $T$ and $T'$ 
cancel in the sum.
We see that we can erase from the sum all trees in which the two border vertices have the same color, but different degrees. 

{\bf B. Symmetrizing an interior stretch.} Choose a stretch 
formed by an even number of consecutive interior white vertices 
in the real part of~$T$. 
Pick the two middle vertices of the stretch. If their degrees are different 
we can interchange these two vertices together with the forests growing on them. 
We obtain a tree~$T'$ such that $\eps(T') = -\eps(T)$. Thus, $T$ and $T'$ cancel in the sum, 
so we can disregard all trees having the same length of the real part as $T$
and in which the two middle white vertices of the stretch have different degrees. 
Now assume that the degrees of the two middle white vertices of the stretch are equal. 
Then, we look at the two white vertices surrounding these two and perform the same operation. 
By the same argument as above, we can disregard all trees~$T$ where the degrees of the 4 white vertices are not symmetric. 
Continuing in the same way we see that, among the trees with a given length of the real part, 
we can erase from the sum all the trees in which the degrees of the real white vertices in a chosen stretch are not symmetric. 
Of course, the same considerations apply to black vertices as well. 

Now consider four cases.

{\bf 1.} If $T$ has only one real vertex we just leave~$T$ in the sum. In this case $\eps(T) = \omega(T) = 1$.

{\bf 2.} The tree $T$ has more than one real vertex and 
the border vertices of $T$ are both white. 

The first possibility is that there is an even number of white interior vertices and an odd number of black interior vertices. By (A) we can assume that the degrees of the border vertices are equal. By (B) we can assume that the degrees of the white vertices are symmetric. Also by (B) we can exclude the leftmost black vertex and assume that the other black vertices have symmetric degrees. The signs of the remaining trees are all equal to $1$. Indeed, because of the symmetry,
all disorders, both black and white, come in pairs. 
Since we are trying to symmetrize the real part sequence of the tree we 
now replace $T$ 
by a new tree $T'$ by moving the leftmost black vertex to the middle position in the real part. 
This 
changes the sign of the tree, but we must still count it 
with the original sign $\omega(T') = \eps(T)$. 

The second possibility is that there is an even number of black interior vertices and an odd number of white interior vertices. By (A) we can assume that the degrees of the border vertices are equal. By (B) we can assume that the degrees of the black vertices are symmetric. Also by (B) we can exclude the leftmost interior white vertex and assume that the other interior white vertices have symmetric degrees. The signs of the remaining trees are all equal to $-1$. Indeed, because of the symmetry, all disorders, 
both black and white, come in pairs, except for the disorders involving a border vertex and the excluded interior white vertex. 
Since the degrees of the border vertices are equal, there is exactly one exceptional disorder like that, 
so the total number of disorders is odd. 
As before, to symmetrize the real part sequence of the tree we replace $T$ by a new tree $T'$ by moving the leftmost interior white vertex to the middle position in the real part. This 
changes the sign of the tree, but we must still count it with the original sign $\omega(T') = \eps(T)$. 

Thus, we see that in both situations all trees with non-symmetric real part sequences cancel out and we are only left with trees 
having symmetric real part sequences. Every tree like that is counted with weight~$\omega$. 

{\bf 3.} The tree $T$ has more than one real vertex and its border vertices are both black. This case is treated in the same way
as the previous one. 

{\bf 4.} One border vertex of~$T$ is white and one is black. In this case there are as many white vertices as black ones. 

If the number of interior white (and therefore black) vertices is even, we can assume by~(B) that their degrees are symmetric. 
Consider the tree $T'$ obtained from~$T$ by a rotation by $180^\circ$. We have $\eps(T) = \eps(T')$. Indeed, as we revert the orientation of the real line all disorders appear or disappear in pairs because of the symmetry. On the other hand, among the trees~$T$ and $T'$ one is a white side tree while the other is a black side tree. Thus, these trees 
cancel in the sum. 

If the number of interior white (and therefore black) vertices is odd, we can exclude the interior vertices closest to the border vertices of the same color and assume by~(B) that the degrees of the remaining interior vertices are symmetric. Consider the tree $T'$ obtained from~$T$ by a rotation by $180^\circ$. Once again, we have $\eps(T) = \eps(T')$. Indeed, as we revert the orientation of the real line all disorders appear or disappear in pairs because of the symmetry, except for the disorders involving the border vertices and their adjacent vertices of the same color. These two exceptional disorders modify the sign by an extra factor of $(-1)^2=1$. Now, as before, among the trees~$T$ and $T'$ one is a white side tree while the other is a black side tree, so they cancel in the sum.

We see that all the trees of the sum have canceled, so the sum is equal to~$0$. Note that in this case there are no trees with symmetric real part sequence.

The lemma is proved. \qed

\paragraph{Proof of Lemma~\ref{Lem:weightzero}.}

If the number $n$ of edges of our trees is odd, every tree has two border vertices of different colors. In that case there are no trees with symmetric real part sequence, so the weights of all trees vanish.

The interesting case is when $n$ is even. In that case we will prove that there is a cancellation between the trees with a unique real vertex and the trees with more than one real vertex.

\begin{figure}[h]
\begin{center}
\includegraphics[width=15em]{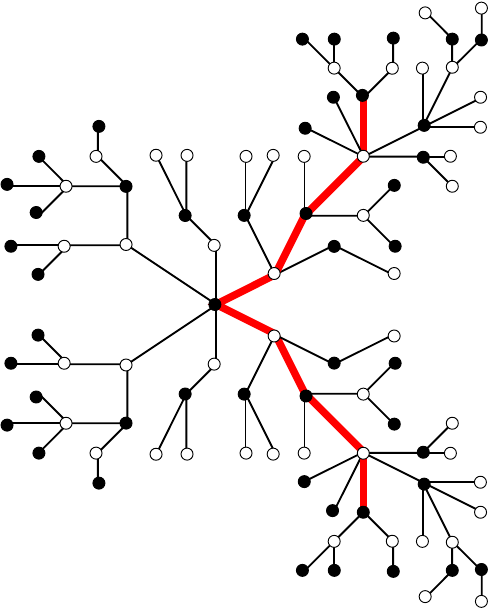}
\caption{A tree with a unique real vertex;
the midlines of its rightmost branches are shown in red.} \label{Fig:midline} 
\end{center}
\end{figure}

Consider a tree~$T$ with a unique real vertex. Take the two symmetric trees growing on it closest to the positive direction of the real axis. In each of these trees construct the {\em midline}: the sequence of edges that starts at the real vertex and divides into two equal parts the degree of each vertex that it meets, see Figure~\ref{Fig:midline}. 

\begin{figure}[h]
\begin{center}
\includegraphics[width=\textwidth]{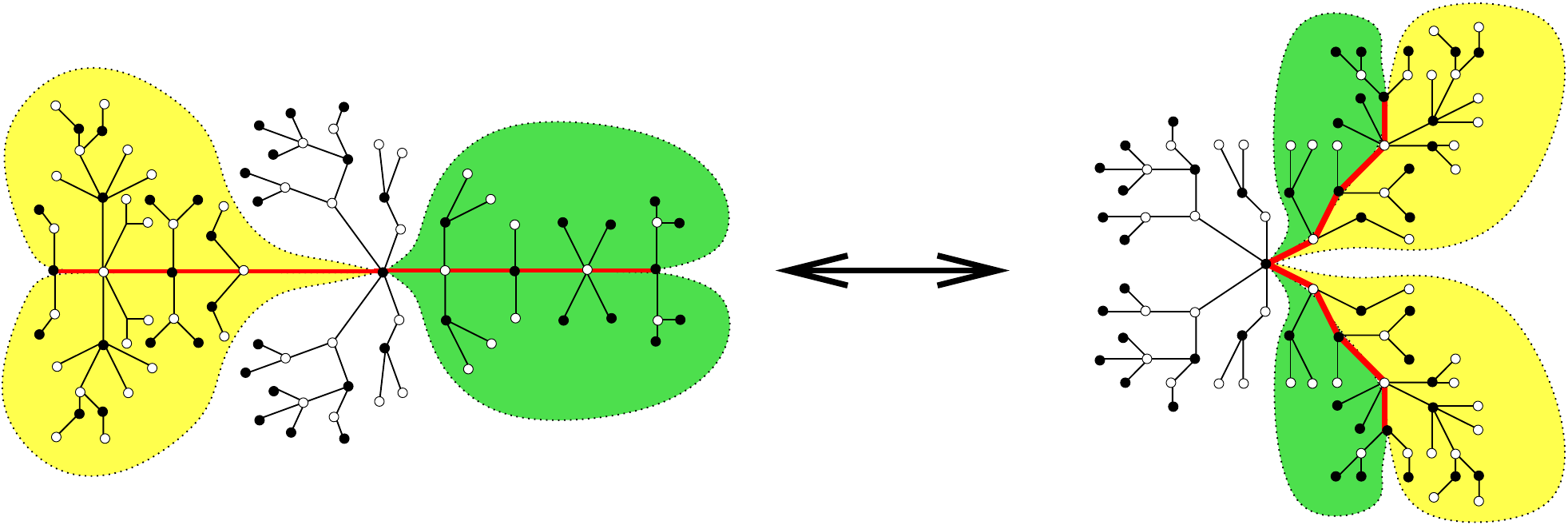}
\caption{A bijection between trees with a unique real vertex and trees with a symmetric real part sequence composed of more than one vertex.} \label{Fig:cutandpaste}
\end{center}
\end{figure}

Denote by $A$ and $A'$ the half-trees closest to the real axis and by $B$ and $B'$ the half-trees separated from the axis by the midlines. Now assemble $A$ and $A'$ into a new real tree along the positive direction of the real axis; the midline will follow the real axis. Similarly, assemble $B$ and $B'$ into a new real tree along the negative direction of the real axis; the midline will follow the real axis. We have obtained a new tree~$T'$ with more than one real vertex and symmetric real part sequence,
see Figure~\ref{Fig:cutandpaste}. 
This construction establishes a bijection between the set of trees with a unique real vertex and 
the set of trees with a symmetric real part sequence composed of more than one vertex
(in both cases we consider only trees with $n$ edges).

If, among the trees $T$ and $T'$, one is a white side tree and the other a black side tree, then $\omega(T) = \omega(T')=1$, so the trees cancel in the sum. If both trees are white side or both are black side, then $\omega(T)=1$ while $\omega(T')=-1$, so the trees cancel as well.

The lemma is proved. \qed

\section{Proof of the invariance theorem}\label{proof}

Let $k$ and $n$ be two positive integers, $k < n$, and let $(\Lambda_1, \dots, \Lambda_k)$ 
be a sequence of partitions of~$n$ such that 
$$
\sum_{i=1}^k (n - l(\Lambda_i)) = n-1.
$$
In view of the equality between the $s$-numbers of real normalized polynomials 
and increasing real polynomial dessins (Corollary~\ref{signed_equality}), 
the invariance theorem (Theorem~\ref{Thm:invariance}) is equivalent to the following proposition. 

\begin{proposition}\label{key-proposition} 
For any integer $1 \leq i \leq k - 1$, the $s$-number
of increasing 
real polynomial dessins of degree $n$
and type 
$$(\Lambda_1, \ldots, \Lambda_{i - 1}, \Lambda_i, \Lambda_{i + 1}, \Lambda_{i + 2}, \ldots, \Lambda_k)$$
is equal to the $s$-number 
of increasing 
real polynomial dessins of degree $n$ 
and type 
$$(\Lambda_1, \ldots, \Lambda_{i - 1}, \Lambda_{i + 1}, \Lambda_i, \Lambda_{i + 2}, \ldots, \Lambda_k).$$ 
\end{proposition} 

The rest of the section is dedicated to a proof of this proposition. The proof uses the invariance theorem for black and white trees (Theorem~\ref{tree-invariance}). To match the notation, we 
put $\Lambda_i = \Lambda_b$ and $\Lambda_{i+1} = \Lambda_w$. 

Let $w_1 < \ldots < w_k$ be real numbers, 
and consider 
be the set of (homeomorphism classes of) increasing real polynomial
dessins $\Gamma$ of degree $n$ and type
$$
(\Lambda_1, \ldots, \Lambda_{i - 1}, \Lambda_b, \Lambda_w, \Lambda_{i + 2}, \ldots, \Lambda_k).
$$

For each real polynomial dessin $\Gamma$
we color in black the vertices labelled with $i$ (they correspond to $\Lambda_b$), and in white the vertices labelled with $i + 1$ (they correspond to $\Lambda_w$).


These black and white vertices together with the edges of type $i \to i +1$ of $\Gamma$ 
form a collection of black and white trees embedded into~$\C$. Some of them are real; others split into pairs of trees that are complex conjugate (symmetric to each other). To distinguish between the two cases we use, 
until the end of this section, the terms {\em real trees} 
for real black and white trees of Definition~\ref{Def:tree} and {\em 
imaginary trees}
for black and white trees 
without vertices on the real axis. 

Given a real polynomial dessin $\Gamma$ as above, 
we can contract all its $i \to i+1$ edges in order to obtain a new real polynomial dessin $\Gh$. Thus,
each real and 
imaginary 
tree gets contracted to a vertex. 
The vertices obtained from contracted trees are labeled with $i$. 
The vertices that were labelled with indices $j > i + 1$ are now labeled with $j - 1$. 
Moreover, we enhance each vertex of $\Gh$ labeled by $i$ with the couple $(n_b, n_w)$ 
of partitions $n_b$ and $n_w$ describing the degrees of white and black vertices of the real or 
imaginary tree contracted to this vertex. 

\begin{definition}\label{Def:new_enhanced} 
An {\em $i$-enhanced dessin} is an increasing real polynomial dessin with an extra data, for each vertex $v$ labelled by $i$, 
of a pair $(n_b, n_w)$ of partitions of $(\deg v)/2$, where $\deg v$ is the degree of $v$.
The vertices labelled with $i$ of an $i$-enchanced dessin are called {\em special}. 
\end{definition} 

\begin{definition} \label{Def:EnhancedDisorder}
Let $\Gh$ be an $i$-enhanced dessin.
A {\em disorder} for vertices with label $j \ne i$ is defined in the usual way. For special vertices, however, we have a special definition of disorders. Consider two real special vertices $v_1<v_2$. 
Choose
an element $n_1$ of one of the partitions
assigned to $v_1$ and an element $n_2$ of the partition of the 
same color assigned to $v_2$. The couple $(n_1, n_2)$ is called a {\em special disorder} if $n_1 > n_2$.
\end{definition} 

\begin{definition} \label{Def:EnhancedSign}
The {\em sign of an enhanced dessin} $\Gh$ is $\eps(\Gamma) = (-1)^d$, where $d$ is the total number of disorders of $\Gh$, both special and ordinary. 
\end{definition}

Denote by $\daleth$ the collection of (homeomorphism classes of) $i$-enhanced dessins 
of type $\Lambda_1, \dots, \Lambda_{i-1}, \Lambda, \Lambda_{i+2}, \dots , \Lambda_k$ 
such that if we take the union over the special vertices of the black partitions $n_b$ 
we get $\Lambda_b$ and if we take the union of all white partitions $n_w$ we get $\Lambda_w$. 
Thus, any $i$-enhanced dessin obtained by contraction from 
an increasing real polynomial dessin of type $(\Lambda_1, \dots, \Lambda_k)$ belongs to $\daleth$. 

\begin{lemma} \label{Lem:EnhancedSign}
If an $i$-enhanced dessin $\Gh \in \daleth$ is obtained from a dessin $\Gamma$ by contraction  then
$$
\eps(\Gamma) = \eps(\Gh) \times \prod_T \eps(T),
$$
where the product is taken over the contracted real trees in $\Gamma$.
\end{lemma}

\paragraph{Proof.} 
Let $d(\Gh)$ be the number of disorders in $\Gh$. Let $\sum_T d(T)$ be  the sum of the numbers of disorders of all real trees in $\Gamma$ that got contracted. Finally, let $D$ be the number of pairs of vertices $(v_1, v_2)$ in two distinct contracted trees $T_1$ and $T_2$ such that 
\begin{itemize}
\item the tree $T_1$ lies to the left of $T_2$,
\item the degree of $v_1$ is greater than the degree of $v_2$,
\item either $v_1$ or $v_2$ (or both) is not real.
\end{itemize}
Then the number of disorders of $\Gamma$ is equal to 
$$
d(\Gh) + \sum_T d(T) - D.
$$
It is easy to see that $D$ is always even, since the pairs $(v_1,v_2)$ come in couples of quadruples of complex conjugate vertices. Thus
$$
\eps(\Gamma) = \eps(\Gh) \times \prod_T \eps(T)
$$
as claimed.
\qed

Now let us study the set of increasing real polynomial dessins $\Gamma$ that contract to a given $\Gh \in \daleth$. 
Let $V_\R$ be the set of special real vertices of $\Gh$, and let $V_+$ be the set of its special vertices 
with a positive imaginary part. 
For convenience, 
for every special vertex $v \in V_\R \cup V_+$, 
we mark an adjacent edge 
in the following way. 
For any vertex $v \in V_\R$, we mark the edge to the right of the vertex on the real axis. 
For any vertex $v \in V_+$, we mark at random an edge of type $i \to i+1$ (or $i \to \infty$ 
if $i$ is the largest critical value).

Assign a set $T_v$ to each special vertex of~$\Gh$ as follows.
\begin{itemize}
\item
If $v$ is in $V_+$ and the corresponding partitions are $(n_b, n_w)$, we assign to $v$ the set $T_v$ of 
imaginary trees with one marked white half-edge ({\it i.e.}, a half-edge adjacent to a white vertex), such that the degrees of black and white vertices of the tree are given by $n_b$ and~$n_w$. 

\item
If $v \in V_\R$ is a real vertex, $(n_b, n_w)$ the corresponding partitions, and the marked edge to the right of $v$ is of type $i \to i+1$ or $i \to \infty$, then we assign to $v$ the set $T_v$ of {\em white side} real trees with degrees of black and white vertices given by $n_b$ and $n_w$. 

\item
Finally, if $v \in V_\R$ is a real vertex, $(n_b, n_w)$ the corresponding partitions, and the marked edge to the right of $v$ is of type $i-1 \to i$ or $\infty \to i$, then we assign to $v$ the set $T_v$ of {\em black side} real trees with degrees of black and white vertices given by $n_b$ and $n_w$. 
\end{itemize}

\begin{lemma} \label{Lem:set}
The set of increasing real polynomial dessins $\Gamma$ of type $(\Lambda_1, \dots, \Lambda_k)$ that contract to $\Gh$ is in a one-to-one correspondance with the product 
$$
\prod_{v \in V_\R \cup V_+} T_v.
$$
\end{lemma}

\paragraph{Proof.} The dessin $\Gamma$ is constructed by inserting into each special vertex $v$ the corresponding tree. The real trees are inserted in the natural way: the marked edge of $\Gh$ is glued to the rightmost vertex of the real tree, and this vertex has the appropriate color by the construction of our sets of trees. For a vertex $v \in V_+$, we insert the tree in such a way that the marked edge of $\Gh$ is glued to the white vertex of the marked white half-edge of the plane tree, right after this half-edge in the counterclockwise direction. For a vertex $v$ with negative imaginary part, the trees are glued so that the invariance of the dessin by complex conjugation is preserved. \qed

Let $\Gh \in \daleth$ be an $i$-enhanced dessin. 
For a vertex $v \in V_+$, denote by $m_v$ the number of imaginary trees with partitions $(n_b, n_w)$ 
given by the vertex and with one marked white half-edge.
For a vertex $v \in V_\R$, denote by $m_v$ the $s$-number of real black side trees 
with partitions $(n_b, n_w)$ given by the vertex~$v$. 

\begin{remark}
Note that, by Theorem~\ref{tree-invariance}, the $s$-number of real {\em black} side trees with partitions $(n_b, n_w)$ 
is the same as the $s$-number of real {\em white} side trees with partitions $(n_b, n_w)$. 
This equality makes it possible to define the number $m_v$ in an invariant way 
and is the main reason why the proof of invariance goes through. 
\end{remark} 

\begin{lemma} \label{Lem:snumber}
The s-number of increasing real polynomial dessins $\Gamma$ of type $(\Lambda_1, \dots, \Lambda_k)$
that contract to $\Gh$ is equal to 
$$
\eps(\Gh) \prod_{v \in V_\R \cup V_+} m_v.
$$
\end{lemma}

\paragraph{Proof.} The set of dessins $\Gamma$ was described in Lemma~\ref{Lem:set}. The sign of each dessin like that is given in Lemma~\ref{Lem:EnhancedSign}. The lemma follows by combining the two statements. \qed

\paragraph{Proof of Proposition~\ref{key-proposition}.}
The $s$-number of increasing real polynomial dessins $\Gamma$ 
of type $(\Lambda_1, \dots, \Lambda_{i-1}, \Lambda_b, \Lambda_w, \Lambda_{i+2}, \dots, \Lambda_k)$  
that contract to $\Gh$ is given by Lemma~\ref{Lem:snumber}. It is clear that this number does not depend on the order of the partitions $\Lambda_b$ and $\Lambda_w$. 
More precisely, if, at each vertex of $\Gh$, we replace the pair of partitions $(n_b, n_w)$ by $(n_w, n_b)$, 
the $s$-number of Lemma~\ref{Lem:snumber} 
stays the same. 

Since the equality holds for each $i$-enhanced dessin, 
we conclude that the $s$-number of 
increasing real polynomial dessins of type 
$$
(\Lambda_1, \dots, \Lambda_{i-1}, \Lambda_b, \Lambda_w, \Lambda_{i+2}, \dots, \Lambda_k)
$$ 
is the same as for the type 
$$
(\Lambda_1, \dots, \Lambda_{i-1}, \Lambda_w, \Lambda_b, \Lambda_{i+2}, \dots, \Lambda_k). 
$$ 
\vspace{-3em} 

\qed 

\vspace{-1em}

\section{Generating series}\label{generating}
In this section we study the generating series for 
$s$-numbers of real normalized polynomials 
and prove Theorems~\ref{Thm:fg}, \ref{Thm:evenvanishing}, \ref{Thm:oddvanishing},  and~\ref{Thm:logasymp}.

Let $\lambda_1, \dots, \lambda_k$ be a fixed list of partitions. We are interested in the $s$-number 
of real normalized polynomials with $k$ branch points $w_1, \dots, w_k$ 
having reduced ramification types $\lambda_1, \dots, \lambda_k$, respectively, 
and with $m$ more simple branch points $w_{k+1}, \dots, w_{k+m}$. 
Since the $s$-number does not depend on the multiplicity of the branch points on the real line, 
we 
 place the branch points so that $w_1 > w_2 > \dots > w_{k+m}$. 
 
Let~$P$ be a real normalized polynomial whose branch points are all real, 
and consider its affine dessin, that is, the real polynomial dessin $\Gamma_P$ without the $\infty$-vertex. 

Choose a real number $\alpha$ lying between $w_k$ and $w_{k+1}$. To prove the theorems,
we 
extract from the affine dessin several parts that contain all the required information about the polynomial. 
\begin{itemize}
\item
The preimage of the interval $(\alpha, \infty)$. In general this preimage consists of several connected components. We will call a connected component {\em interesting} if it contains either a real point or a critical point of~$P$. The union of interesting components is called the {\em base} of the affine dessin. 
\item 
The real part of the preimage of the interval $(-\infty, w_{k+1}]$. The real part of $P^{-1} \Bigl( (-\infty, w_{k+1}] \Bigr)$ also may consist of several connected components. We call them {\em chains}. 
The real part of $P^{-1} \Bigl( (-\infty, w_{k + 1}] \Bigr)$ 
automatically contains all the critical preimages of $w_{k+1}, \dots, w_{k+m}$. 
\end{itemize}

\begin{figure}[h]
\begin{center}
\ 
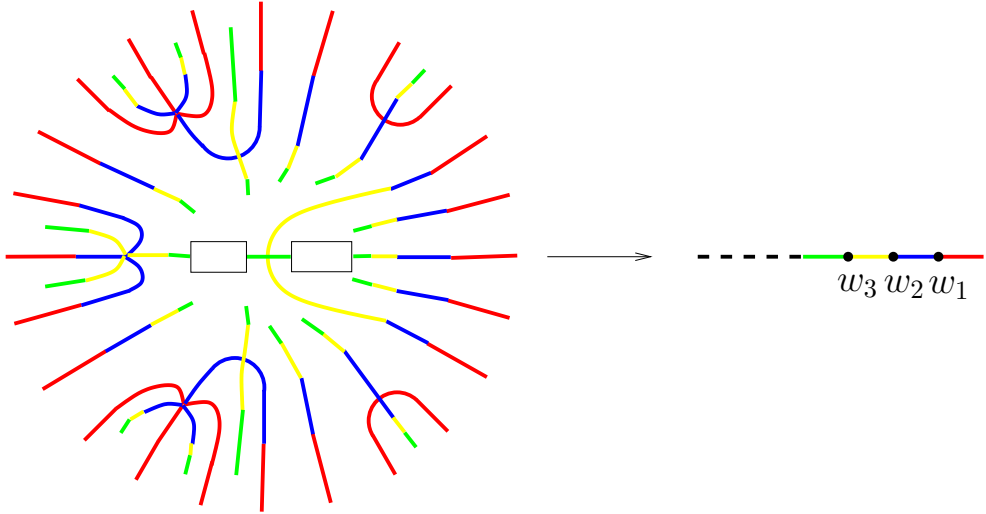

\caption{ The complete preimage of $(\alpha, \infty)$.  Here $n=24$, $k=3$, $\lambda_1 = (2,2,1,1)$, $\lambda_2 = (2,1,1)$, $\lambda_3 = (1)$. The black boxes contain the chains.} 
\label{Fig:preimage}
\end{center}
\end{figure}

\begin{figure}[h]
\begin{center}
\ 
\includegraphics[width=20em]{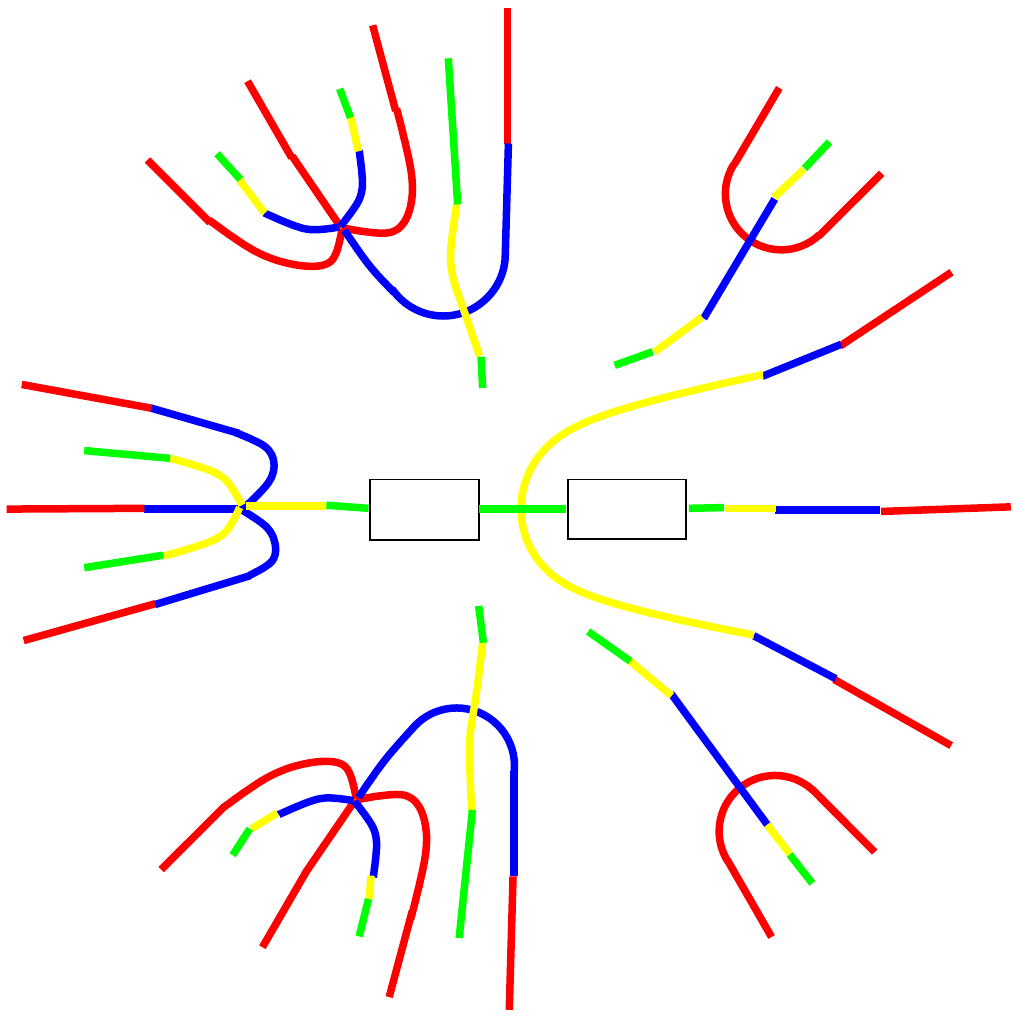}

\caption{The base of the affine dessin. The black boxes contain the chains.} \label{Fig:base}
\end{center}
\end{figure}

The full preimage of $(\alpha, \infty)$ and the base extracted from it are shown in Figures~\ref{Fig:preimage} and~\ref{Fig:base} respectively. The green lines in the picture are the preimages of $(\alpha, w_k)$.
Some of them are contained in the regions between the red rays. 
These green lines 
lead directly to the $n$ preimages of $-\infty$ in the affine dessin without encountering any critical points. 
Other green lines lie in the regions that contain a black box. Those 
are connected to a chain. 

\begin{definition}
Suppose that a base of the affine dessin and its chains are given.
We 
call {\em regions} the connected components into which the base divides the upper half-plane $\{ \Re x > 0 \}$. 
For each chain $c$, let $p_c$ be the number of preimages of $\alpha$ in the region adjacent to~$c$. 
Then, on each chain we can single out the $p_c$ maxima to which these preimages are connected. 
A chain with $p_c$ maxima singled out 
is called {\em marked}.
\end{definition} 

Theorems~\ref{Thm:fg}, \ref{Thm:evenvanishing}, and \ref{Thm:oddvanishing}
are proved by counting the $s$-number of real normalized polynomials 
with a given base and then performing a summation over the finite set of possible bases. 
Before proceeding with the proofs we state three simple lemmas. 

\begin{lemma} \label{Lem:BaseToTree}
An affine dessin can be uniquely recovered if its base and marked chains are given. 
\end{lemma}

\paragraph{Proof.} Once we have connected the preimages of $\alpha$ to the chains according to the markings, all the vertices of the affine dessin of degree greater than~2 and the edges between them are drawn. The remaining part of the affine dessin is a union of rays consisting of vertices of degree~2 and ending at vertices labeled $\pm \infty$. The first interval of each ray like that is uniquely determined by the vertex structure of the affine dessin, namely, by the fact that the edges $(i \rightarrow i + 1)$ and $(i - 1 \rightarrow i)$
alternate at each vertex of color~$k$. 
Once we know the first interval of a ray, it can be extended in a unique way. \qed

\begin{lemma} \label{Lem:D}
Introduce the differential operator $D = q \frac{d}{dq}$. Then, we have 
\begin{align*}
Dq & = q, \\
Df & = q(1-f^2),\\
Dg & = -qfg.
\end{align*}
\end{lemma}

\paragraph{Proof.} This is just a computation. \qed

\begin{notation}
Put 
\begin{align*}
f_p &= \frac1{2^p \, p!} (D-1) (D-3) \cdots (D-2p+1) \, f, \\
g_p &= \frac1{2^{p} \, p!} D(D-2) \cdots (D-2p+2) \, g.
\end{align*}
\end{notation}

\begin{lemma} \label{Lem:fpgp}
The series $f_p$ {\rm (}respectively, $g_p${\rm )} is the generating series for the numbers of alternating permutations 
of odd {\rm (}respectively, even{\rm )} length with $p$ distinguished maxima.
\end{lemma} 

\paragraph{Proof.} The operator~$D$ multiplies the coefficient of $q^m$ in a generating series by~$m$. An alternating permutation of odd length~$m$ has $(m-1)/2$ maxima. Therefore, to choose~$p$ maxima we need to multiply the number of alternating permutations by 
$$
\frac{\frac{m-1}2 \cdot \frac{m-3}2 \cdot \dots \cdot \frac{m-2p+1}2}{p!}.
$$
An alternating permutation of even length~$m$ has $m/2$ maxima. Therefore,
to choose~$p$ maxima we need to multiply the number of alternating permutations by 
$$
\frac{\frac{m}2 \cdot \frac{m-2}2 \cdot \dots \cdot \frac{m-2p+2}2}{p!}.
$$
\qed

\paragraph{Proof of Theorem~\ref{Thm:fg}.}
Let $\lambda_1, \dots, \lambda_k$ be a given list of partitions. Consider a possible base of the affine dessin $P^{-1}(\R)$. Denote by $\eps$ the sign of the base, that is, $(-1)$ to the number of disorders present in the base. Suppose that the base leaves $b$ gaps for chains and that there are $p_1, p_2, \dots, p_b$ maxima to choose in the $b$ chains to be connected with the preimages of~$\alpha$ in the base. 

\begin{figure}[h]
\begin{center}
\ 
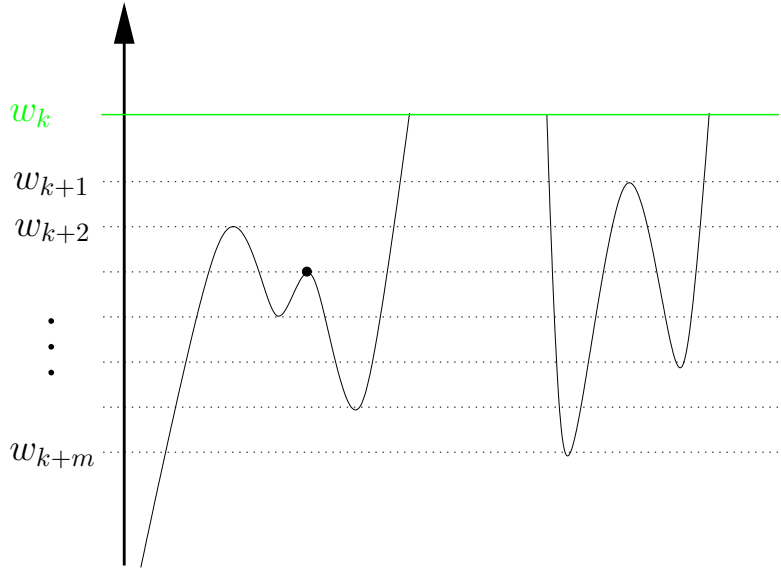

\caption{The chains. Here $n$ is odd, so that the leftmost chain extends to~$-\infty$. There are $m=7$ simple branch points divided between $b=2$ chains. The numbers of selected maxima is $p_1 = 1$ and $p_2 = 0$.} 
\end{center}
\end{figure} 

We claim that the $s$-number of polynomials with the given base and $m$ 
additional simple branch points equals the coefficient of $q^m/m!$ in the power series 
\begin{align}
\eps \cdot \prod_{i=1}^b f_{p_i}  & \qquad \mbox{ for } n \mbox{ even,} \label{Eq:even}\\
\eps \cdot g_{p_1} \prod_{i=2}^b f_{p_i} & \qquad \mbox{ for } n \mbox{ odd}. \label{Eq:odd}
\end{align}
Indeed, it follows from Lemma~\ref{Lem:fpgp} that the 
products $\prod_{i=1}^b f_{p_i}$ and $g_{p_1} \prod_{i=2}^b f_{p_i}$ correctly count the chains with distinguished maxima, taking into account the disorders within each chain. The sign $\eps$ also takes into account the disorders in the base. Note that all the real preimages of critical values are contained in the union of the base and the chains. Thus, every disorder takes place either 
within the base (for critical values between $w_1$ and $w_k$), or within a chain, or between two chains. 
The disorders with each chain are accounted for by the signs in the series $\prod_{i=1}^b f_{p_i}$ and $g_{p_1} \prod_{i=2}^b f_{p_i}$. 
The number of disorders between two chains is always even. Indeed, a disorder between two chains appears between 
a simple critical point on some level~$w_{k+i}$ in the leftmost chain and a simple preimage of~$w_{k+i}$ in the rightmost chain. 
But the rightmost chain crosses the level $w_{k+i}$ an even number of times, so there is an even number of disorders. 

Now, Lemma~\ref{Lem:D} implies that~\eqref{Eq:even} is a polynomial in $q$ and $f$, while~\eqref{Eq:odd} is $g$ times a polynomial in $q$ and $f$. The generating series 
$\Fe_{\lambda_1, \dots, \lambda_k}$ (respectively, $\Fo_{\lambda_1, \dots, \lambda_k}$) is obtained as a sum of expressions~\eqref{Eq:even} (respectively,~\eqref{Eq:odd}) over the finite set of all possible bases of the affine dessin. 
Thus, $
\Fe_{\lambda_1, \dots, \lambda_k}$ is a polynomial in $q$ and $f$, 
while $
\Fo_{\lambda_1, \dots, \lambda_k}$ is $g$ times a polynomial in $q$ and $f$.
\qed 

\paragraph{Proof of the {\em if} part of Theorem~\ref{Thm:evenvanishing} (nonvanishing for even degree polynomials).}
Consider $k$ partitions $\lambda_1, \dots, \lambda_k$ and assume that in each of them every even number appears an even number of times and at most one odd number appears an odd number of times. Let $s$ be the number of partitions in which one of the numbers appears an odd number of times. 

Given a partition $\lambda$, denote by $[\lambda/2]$ the partition defined by the following rule: if an integer $a$ appears $n_a$ times in $\lambda$ then it appears $[n_a/2]$ times in $[\lambda/2]$, where $[ \cdot]$ is the integer part. For instance, if $\lambda = (1,1,1,1,1,2,2,3,4,4,6,8)$ then $[\lambda/2] = (1,1,2,4)$.

Denote by~$\ell$ the total number of elements in the partitions $[\lambda_i/2]$:
$$
\ell = l([\lambda_1/2]) + \cdots + l([\lambda_k/2]).
$$

\begin{figure}[h]
\begin{center}
\ 
\includegraphics[width=25em]{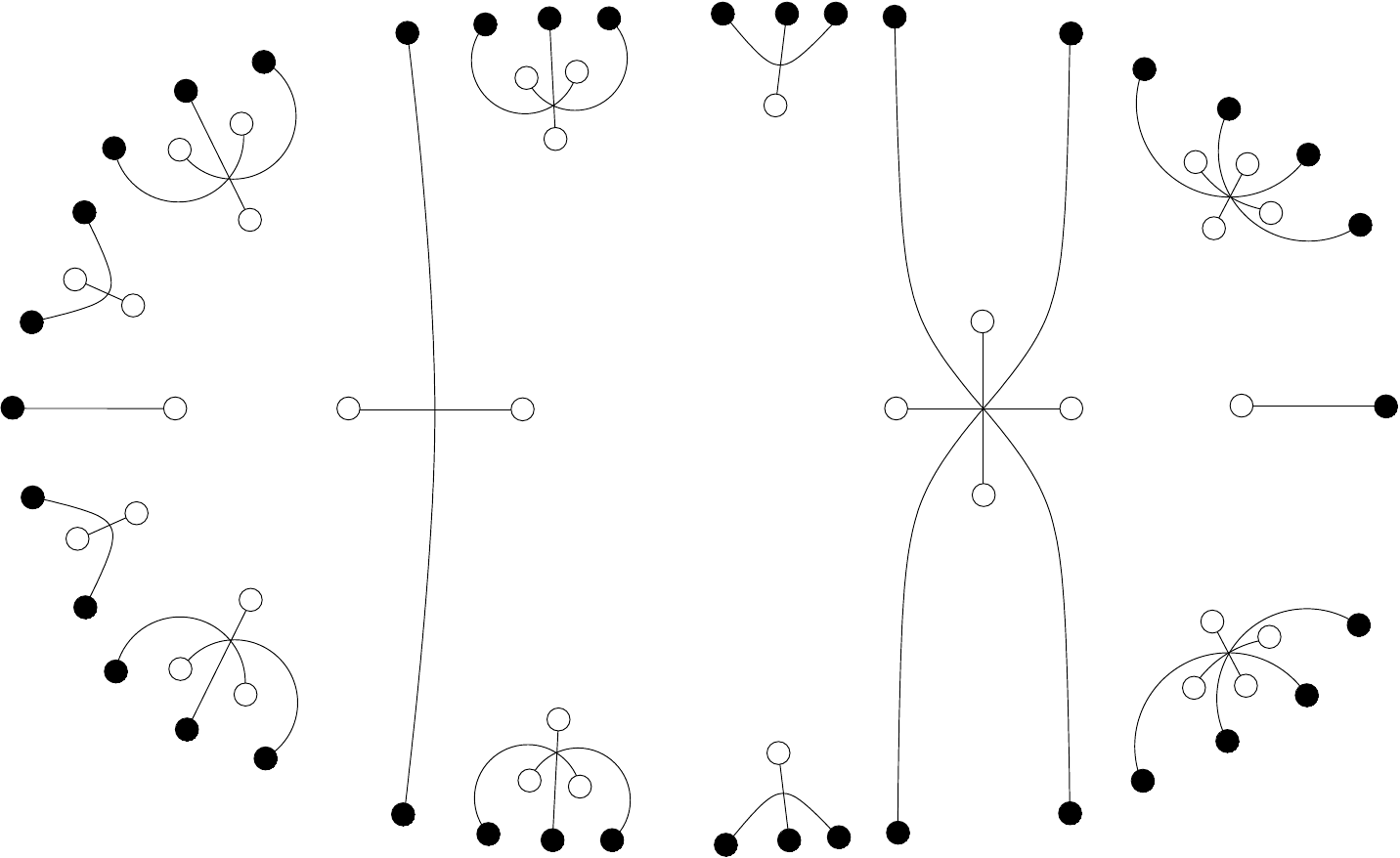}
\end{center}
\caption{A simple base for even~$n$. The preimages of $+\infty$ are represented as black dots; those of $\alpha$ as white dots. In this example $s=2$, $\ell=5$.} \label{Fig:SimpleBaseEven}
\end{figure}

Let us call a base {\em simple} if all critical points lie on different connected components of the base and if there is at most one real critical point on each level $w_1, \dots, w_k$, see Figure~\ref{Fig:SimpleBaseEven}. In particular, a simple base contains $s$ real critical points and $\ell$ pairs of complex conjugate critical points.

We know that the series $\Fe_{\lambda_1, \dots, \lambda_k}$ is a polynomial in $q$ and~$f$. 
We claim that {\em each simple base gives rise to a nonzero contribution to the coefficient of $q^\ell f^{\ell+s+1}$ 
in this polynomial; this contribution has sign $(-1)^\ell$; moreover, no other base contributes to this coefficient.} 
To prove this claim, we, 
first, study the power of~$q$, then the power of~$f$, and then the sign. 

{\bf The power of~$q$.}
The series $\Fe_{\lambda_1, \dots, \lambda_k}$ is a sum of expressions of the form
$$
\eps \cdot \prod_{i=1}^b f_{p_i},
$$
where $p_i$ is the number of distinguished maxima in the $i$th chain. Recall that $Df = - q (f^2-1)$, $Dq = q$. The highest possible power of~$q$ in 
the above expression equals $\sum p_i$, because every time we apply the operator~$D$ 
we increase the power of~$q$ by at most~1. 
Thus, the highest possible coefficient of~$q$ 
is obtained from the bases with the largest possible number of distinguished maxima. 
Moving along the affine dessin from a distinguished maximum of a chain towards a preimage of $\alpha$ and beyond, 
we 
eventually arrive at a non-real critical point. Thus, the greatest possible number of distinguished maxima is achieved 
when as few critical points as possible are real and when there is a distinguished maximum assigned to each non-real critical point. 
The number of distinguished maxima 
is then 
equal to $\ell$, so this is also the highest possible degree in~$q$. 

{\bf The power of~$f$.} From now on we 
restrict ourselves to bases that contribute to the coefficient of $q^\ell$. 
We know that the contribution of a base like that is obtained by applying $\ell$ copies of~$D$ to $b$ copies of~$f$, 
where $b$ is the number of chains. 
Applying the operator~$D$ increases the power of $f$ by at most~1. Thus, we 
get the highest possible power of~$f$ if we start with as many chains as possible. 
Since we only have $s$ real critical points available on levels $w_1$ to $w_k$, 
there 
are 
at most $s+1$ chains. In that case, the largest power of~$f$ 
is equal to $\ell+s+1$ and the base 
is simple. 

{\bf The sign.} 
In a simple base there is exactly one real critical point between two chains. Therefore, 
this critical point is a local maximum. It follows that there is an even number of simple preimages lying to the right of this critical point on the same level. Thus, the number of disorders in a simple base is even, so its sign equals $\eps= 1$. 
The coefficient of the monomial $q^\ell f^{\ell+s+1}$ in $\prod_{i=1}^{s+1} f_{p_i}$ is equal to $(-1/2)^\ell$. 
Indeed, we have a change of sign and a division by~2 every time we apply the operator~$D$, and there are $\ell$ operators to apply. 
Thus, we see that every simple base contributes exactly $(-1/2)^\ell$ to the coefficient of $q^\ell f^{\ell+s+1}$. 

Simple bases obviously exist; thus we see that the coefficient of 
$q^\ell f^{\ell+s+1}$ is nonzero. For completeness, let us compute this coefficient. The number of simple bases is given by the number of possible orderings of the pairs of complex conjugate critical points (or, more precisely, their distinguished maxima) and real critical points. One should also take into account the fact that critical points of the same multiplicity lying on the same level are indistinguishable. Denote by $\Aut [\lambda_i/2]$ the number of automorphisms of the partition $[\lambda_i/2]$. Then the coefficient of $q^\ell f^{\ell+s+1}$ equals
$$
\left(-\frac12\right)^\ell  \frac{(\ell+s)!}{\prod\limits_{i=1}^k \Aut [\lambda_i/2]}.
$$

\qed 

\paragraph{Proof of the {\em if} part of Theorem~\ref{Thm:oddvanishing} (nonvanishing for odd degree polynomials).} 
The proof goes along the same lines as above, but with an extra complication. 
We 
skip some details in the parts of the proof strictly analogous to the proof above, but highlight the differences. 

Consider $k$ partitions $\lambda_1, \dots, \lambda_k$ and assume that in each of them at most one even and at most one odd element appears an odd number of times. Let $s$ be the number of partitions in which one of the {\em odd} numbers appears an odd number of times. 

To every partition $\lambda_i$ we assign a sign $\eps_i = \pm 1$ in the following way.
If $\lambda_i$ has an odd element that appears an odd number of times and
no even element that appears an odd number of times, then $\eps_i = -1$.
If $\lambda_i$ has an odd element that appears an odd number of times and an even element that appears 
an odd number of times and if the odd element is greater than the even element then 
$\eps_i = -1$. 
In all other cases we set $\eps_i = 1$. We also let $\eps = \prod_{i=1}^k \eps_i$. 

As before, we introduce the partitions $[\lambda_i/2]$ and denote by $\ell$ their total length
$$
\ell = \sum_{i=1}^k l([\lambda_i/2]).
$$

To avoid confusion we will call {\em crossing} and {\em extremal} the real critical points of even and odd multiplicity respectively. Thus the graph of the polynomial $P$ crosses the horizontal line $y=w_i$ at a crossing critical point, but not an extremal one. 

\begin{figure}[h]
\begin{center}
\ 
\includegraphics[width=30em]{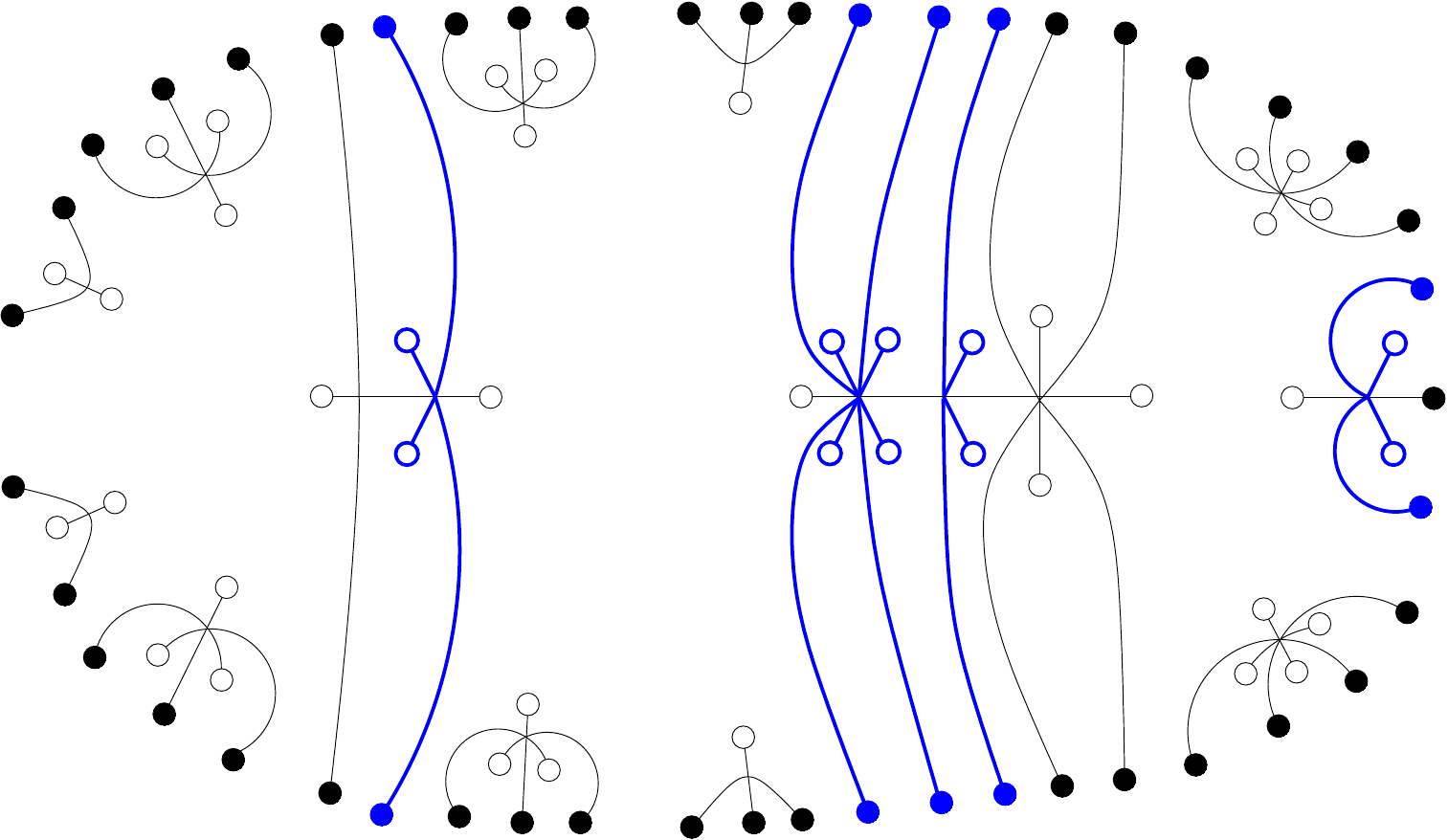}
\end{center}
\caption{A simple base for odd~$n$. The preimages of $+\infty$ are represented as black dots; those of $\alpha$ as white dots. In this example $s=2$, $\ell=5$.} \label{Fig:SimpleBaseOdd}
\end{figure}

Let us call a base {\em simple} if all non-real critical points and all real extremal critical points lie on different connected components of the base and if there is at most one real crossing and at most one real extremal critical point on each level $w_1, \dots, w_k$, see Figure~\ref{Fig:SimpleBaseOdd}. 

A simple base contains $\ell$ pairs of complex conjugate critical points and $s$ real extremal critical points. Each of those critical points is contained in its own connected component of the base. One more connected component of the base contains the real half-line to the right of the rightmost chain. Altogether there are $s+1$ ``real'' and $2 \ell$ ``imaginary'' connected components. In addition to the critical points mentioned above, there is a certain number of real crossing critical points. Those  are distributed among the $s+1$ ``real'' connected components of the base.

We call the {\em skeleton} of a simple base the part obtained by removing 
all the non-real edges adjacent to real crossing critical points. In 
Figure~\ref{Fig:SimpleBaseOdd} 
the edges to be removed are shown in thick blue lines, 
while the skeleton is shown in black.

We know that the series $\Fo_{\lambda_1, \dots, \lambda_k}$ is $g$ times a polynomial in $q$ and~$f$. We claim that {\em each simple base gives rise to a nonzero contribution to the coefficient of $q^\ell f^{\ell+s}$ in this polynomial; the sum of these contributions over all simple bases with the same skeleton has sign $(-1)^\ell \cdot \eps$; moreover, no other base contributes to this coefficient.} 
To prove this claim, we, 
first, study the power of~$q$, then the power of~$f$, and then the sign. 

{\bf The power of~$q$.}
As before, the highest possible coefficient of~$q$ is obtained from the bases with the largest possible number of distinguished maxima, or, in other words, when as few critical points as possible are real and when there is a distinguished maximum assigned to each non-real critical point. The number of distinguished maxima 
is then 
equal to $\ell$, so this is also the highest possible degree in~$q$. 

{\bf The power of~$f$.} From now on we 
restrict ourselves to bases that contribute to the coefficient of $q^\ell$. We know that the contribution of a base like that is obtained by applying $\ell$ times the operator $D$ to one copy of~$g$ (the left-most chain) and several copies of~$f$ (the other chains). Applying the operator~$D$ increases the power of $f$ by at most~1. Thus, we 
get the highest possible power of~$f$ if we start with as many chains as possible. 
Note that two chains have to be separated by at least one local maximum in the base, that is, a real extremal critical point. Since we only have $s$ of those available on levels $w_1$ to $w_k$, there 
are at most $s+1$ chains. In that case, the largest power of~$f$ 
is equal to $\ell+s$ and the base 
is simple. 

{\bf The sign.} Given the skeleton of a simple base, one can ``graft'' the crossing critical points and the blue edges growing out of them in several ways. More precisely, suppose that in the partition $\lambda_i$ an even element appears an odd number of times. In that case (and only in that case) simple bases 
contain a real crossing critical point on level~$w_i$. All real preimages of $w_i$ are contained in the skeleton. One of them may be the extremal critical point (if there is one). At all the other preimages the graph of the polynomial crosses the level $y = w_i$ and therefore there is an odd number of such points. The crossing critical point can be grafted at any of these points. Note that the signs of the bases thus obtained alternate. Indeed, every time we move the crossing critical point one position to the right without jumping over the extremal critical point we destroy exactly one disorder. If we do jump over the extremal critical we destroy one disorder and might create two more, but in both cases the parity of the number of disorders changes. We see that the contributions of the simple bases cancel out, except for the last simple base where the crossing critical point is at the rightmost position. The same reasoning holds for every level~$w_i$. Thus, we conclude that for a given skeleton the contributions of all simple bases cancel out, except for the one simple base in which all crossing critical points are at the rightmost positions on their respective levels, that is, to the right of the last chain. The sign of this simple base equals~$\eps$ (defined at the beginning of the proof).
Indeed, consider a critical level $w_i$ for $1 \leq i \leq k$. It contains at most two real critical points: at most one extremal one that is then necessarily a local maximum, and at most one crossing one. If it has no real critical points, 
it contains no disorders and contributes a sign $+1$. If it has only a crossing critical point, it contains no disorders, 
since the crossing point is in the rightmost position. 
If it has only a local maximum, it contains an odd number of disorders, and therefore contributes a sign $-1$.
If it has both a local maximum and a crossing critical point in the rightmost position,
the parity of disorders depends on which of these two critical points has a greater multiplicity. 

The coefficient of the monomial $g \cdot q^\ell f^{\ell+s}$ in
$$
g_{p_1} \prod_{i=2}^s f_{p_i}
$$
is equal to $(-1/2)^\ell$. It follows that every skeleton 
contributes exactly $(-1/2)^\ell \cdot \eps$ to the coefficient of $q^\ell f^{\ell+s}$. 

For completeness, let us compute this coefficient of  $q^\ell f^{\ell+s}$. It is given by the number of possible orderings of the pairs of complex conjugate critical points (or, more precisely, their distinguished maxima) and real odd critical points. One should also take into account that critical points of the same multiplicity lying on the same level are indistinguishable. Denote by $\Aut [\lambda_i/2]$ the number of automorphisms of the partition $[\lambda_i/2]$. Then, the coefficient of $q^\ell f^{\ell+s}$
equals
$$
\left(-\frac12\right)^\ell \eps \; \frac{(\ell+s)!}{\prod\limits_{i=1}^k \Aut [\lambda_i/2]}.
$$
\qed

\paragraph{Proof of the {\em only if} part of Theorem~\ref{Thm:evenvanishing} (vanishing for even degree polynomials).} 
Assume that in one of the partitions $\lambda_i$ an even element appears an odd number of times. Then, there must be an even real critical point on level~$w_i$. Since the $s$-number of polynomials does not depend on the order of critical values, we may assume that $w_i$ is the global minimum of the polynomial. But the lowest level can only contain odd critical points,
so there are no polynomials at all satisfying the given branching conditions. 

Now assume that in one of the partitions $\lambda_i$ there are two odd elements $a$ and $b$ each of which appears an odd number of times. As before we can assume that $w_i$ is the lowest critical level. We 
divide all dessins under consideration into pairs of dessins with opposite signs. 

Given a dessin $\Gamma$, let us find all of its real vertices corresponding to level $w_i$ and to critical points of multiplicities $a$ and $b$.
There 
is an odd number of vertices of either type. Now make two cuts on the real line to the left and to the right of every chosen vertex.
Remove the chosen vertices from the affine dessin, together with all the edges that grow on them. Now place them back into the affine dessin in the reversed order. We have obtained a new dessin~$\Gamma'$. It is obvious that this operation is an involution: if we reverse the order of the chosen vertices again we 
get back the dessin~$\Gamma$. Note that it is only possible to exchange critical points in that way if they are all  
local minima or all local maxima. In our case, since we have chosen $w_i$ to be the lowest critical level, all critical points 
are local minima. 

Now we claim that the operation $\Gamma \mapsto \Gamma'$ described above changes the parity of the number of disorders. Indeed, consider first the disorders between two chosen vertices. There is an odd number of pairs $(v,w)$, where $v$ is a vertex of multiplicity $a$ and $w$ a vertex of multiplicity~$b$. When we have reversed the order of the vertices, each pair like that has reversed its type: if it was a disorder it is no longer a disorder and if it was not a disorder it has become a disorder. 

Now consider the disorders between a chosen and a non chosen vertex.
Suppose $v$ and $w$ are two chosen vertices that got permuted and $u$ is another vertex with label~$i$. If $u$ does not lie between $v$ and $w$ the number of disorders between $v$ and $w$ on the one hand and $u$ on the other hand does not change. If $u$ lies between $v$ and $w$ it can change by 2 or remain unchanged. Thus the parity of the number of disorders like that has not changed.

To sum up, we see that the total number of disorders has changed parity.
Since we have divided all dessins into pairs of opposite signs, we conclude that the $s$-number of polynomials vanishes. \qed

\paragraph{Proof of the {\em only if} part of Theorem~\ref{Thm:oddvanishing} (vanishing for odd degree polynomials).} 
Assume that in one of the partitions $\lambda_i$ two different even elements appear an odd number of times each. Then, 
there must be at least two even real critical points on level~$w_i$. Since the $s$-number of polynomials does not depend on the order of critical values, we may assume that $w_i$ is the lowest critical level. But the lowest level can only contain one even critical point, so there are no polynomials at all satisfying the given branching conditions.

Now assume that in one of the partitions $\lambda_i$ there are two odd elements $a$ and $b$ each of which appears an odd number of times. The proof repeats literally the proof in the even degree case. \qed 

\paragraph{Proof of Theorem~\ref{Thm:logasymp}.} Given a holomorphic function in the disc $|q| < r$ with a unique singularity on the circle $|q|=r$, it is well-known that the coefficients $a_m$ of its Taylor expansion at 0 satisfy 
$$
\ln |a_m| \sim -m \ln r.
$$
In our case, 
the generating function $F$ is holomorphic on $|q| < \pi/2$ and has exactly two poles at $q = \pm i \pi/2$. 
(Indeed, it follows from the proof of Theorems~\ref{Thm:evenvanishing} and~\ref{Thm:oddvanishing} that $\Fe$ and $\Fo/g$ are polynomials in~$q$ and~$f$ of nonzero degree in~$f$. Both $f$ and~$g$ have poles at $\pm i \pi/2$ and these poles cannot cancel out, because $\Fe$ and $\Fo/g$ are polynomials in $q$ and $f$ with rational coefficients, while $i\pi/2$ is transcendental.) Thus,
we have to apply the property above after dividing the generating function $F$ by~$q$ 
if it is odd and substituting $Q = q^2$. 
We obtain that the logarithmic asymptotic of even (if $F$ is even) or odd (if $F$ is odd) coefficients of $F$ is given by $-m \ln(\pi/2)$. 
Finally, we are actually interested in the logarithmic asymptotic of the coefficients multiplied by $m!$. 
Taking into account that $\ln m! \sim m \ln m$ we see that the factorial ``beats'' the coefficients of $F$ 
so that the logarithmic asymptotic of the $s$-numbers is 
equal to $m \ln m$.
\qed

{\ncsc Institut de Math\'ematiques de Jussieu - Paris Rive Gauche\\[-15.5pt] 

Universit\'e Pierre et Marie Curie\\[-15.5pt] 

4 place Jussieu,
75252 Paris Cedex 5,
France} \\[-15.5pt]

{\it and} {\ncsc D\'epartement de math\'ematiques et applications\\[-15.5pt]

Ecole Normale Sup\'erieure\\[-15.5pt]

45 rue d'Ulm, 75230 Paris Cedex 5, France} \\[-15.5pt]

{\it E-mail address}: {\ntt     ilia.itenberg@imj-prg.fr} 

\vskip10pt

{\ncsc Universit\'e{} de Versailles \\[-15.5pt]

CNRS\\[-15.5pt]

25 avenue des \'Etats-Unis,
78000 Versailles,
France} \\[-15.5pt]

{\it E-mail address}: {\ntt     dimitri.zvonkine@uvsq.fr} 

\end{document}

%% file: disorderP.pdftex_t
\begin{picture}(0,0)%
\includegraphics{disorderP.pdf}%
\end{picture}%
%
%
\setlength{\unitlength}{3947sp}%
\begingroup\makeatletter\ifx\SetFigFont\undefined%
\gdef\SetFigFont#1#2#3#4#5{%
  \reset@font\fontsize{#1}{#2pt}%
  \fontfamily{#3}\fontseries{#4}\fontshape{#5}%
  \selectfont}%
\fi\endgroup%
\begin{picture}(4340,3159)(2792,-2748)
\put(4114,-2149){\makebox(0,0)[lb]{\smash{{\SetFigFont{14}{16.8}{\rmdefault}{\mddefault}{\updefault}{\color[rgb]{0,0,0}$x_1$}%
}}}}
\put(5295,-2136){\makebox(0,0)[lb]{\smash{{\SetFigFont{14}{16.8}{\rmdefault}{\mddefault}{\updefault}{\color[rgb]{0,0,0}$x_2$}%
}}}}
\put(4764,-2142){\makebox(0,0)[lb]{\smash{{\SetFigFont{14}{16.8}{\rmdefault}{\mddefault}{\updefault}{\color[rgb]{0,0,0}$x'_1$}%
}}}}
\put(5670,-2136){\makebox(0,0)[lb]{\smash{{\SetFigFont{14}{16.8}{\rmdefault}{\mddefault}{\updefault}{\color[rgb]{0,0,0}$x'_2$}%
}}}}
\end{picture}%

%% file: 2poly.pdftex_t
\begin{picture}(0,0)%
\includegraphics{2poly.pdf}%
\end{picture}%
%
%
\setlength{\unitlength}{2368sp}%
\begingroup\makeatletter\ifx\SetFigFont\undefined%
\gdef\SetFigFont#1#2#3#4#5{%
  \reset@font\fontsize{#1}{#2pt}%
  \fontfamily{#3}\fontseries{#4}\fontshape{#5}%
  \selectfont}%
\fi\endgroup%
\begin{picture}(9300,4382)(2792,-3499)
\put(3014,-1149){\makebox(0,0)[lb]{\smash{{\SetFigFont{14}{16.8}{\rmdefault}{\mddefault}{\updefault}{\color[rgb]{0,0,0}$w_2$}%
}}}}
\put(3039,-2361){\makebox(0,0)[lb]{\smash{{\SetFigFont{14}{16.8}{\rmdefault}{\mddefault}{\updefault}{\color[rgb]{0,0,0}$w_1$}%
}}}}
\put(8464,-924){\makebox(0,0)[lb]{\smash{{\SetFigFont{14}{16.8}{\rmdefault}{\mddefault}{\updefault}{\color[rgb]{0,0,0}$w_2$}%
}}}}
\put(8489,-2374){\makebox(0,0)[lb]{\smash{{\SetFigFont{14}{16.8}{\rmdefault}{\mddefault}{\updefault}{\color[rgb]{0,0,0}$w_1$}%
}}}}
\put(4576,-3374){\makebox(0,0)[lb]{\smash{{\SetFigFont{14}{16.8}{\rmdefault}{\mddefault}{\updefault}{\color[rgb]{0,0,0}$P_1$}%
}}}}
\put(9801,-3349){\makebox(0,0)[lb]{\smash{{\SetFigFont{14}{16.8}{\rmdefault}{\mddefault}{\updefault}{\color[rgb]{0,0,0}$P_2$}%
}}}}
\end{picture}%

%% file: preimage.pdftex_t
\begin{picture}(0,0)%
\includegraphics{preimage.pdf}%
\end{picture}%
%
%
\setlength{\unitlength}{1579sp}%
\begingroup\makeatletter\ifx\SetFigFont\undefined%
\gdef\SetFigFont#1#2#3#4#5{%
  \reset@font\fontsize{#1}{#2pt}%
  \fontfamily{#3}\fontseries{#4}\fontshape{#5}%
  \selectfont}%
\fi\endgroup%
\begin{picture}(15463,8122)(1127,-9127)
\put(14968,-5628){\makebox(0,0)[lb]{\smash{{\SetFigFont{14}{16.8}{\rmdefault}{\mddefault}{\updefault}{\color[rgb]{0,0,0}$w_2$}%
}}}}
\put(14234,-5611){\makebox(0,0)[lb]{\smash{{\SetFigFont{14}{16.8}{\rmdefault}{\mddefault}{\updefault}{\color[rgb]{0,0,0}$w_3$}%
}}}}
\put(15668,-5644){\makebox(0,0)[lb]{\smash{{\SetFigFont{14}{16.8}{\rmdefault}{\mddefault}{\updefault}{\color[rgb]{0,0,0}$w_1$}%
}}}}
\end{picture}%

%% file: chains.pdftex_t
\begin{picture}(0,0)%
\includegraphics{chains.pdf}%
\end{picture}%
%
%
\setlength{\unitlength}{2368sp}%
\begingroup\makeatletter\ifx\SetFigFont\undefined%
\gdef\SetFigFont#1#2#3#4#5{%
  \reset@font\fontsize{#1}{#2pt}%
  \fontfamily{#3}\fontseries{#4}\fontshape{#5}%
  \selectfont}%
\fi\endgroup%
\begin{picture}(8088,5972)(911,-7698)
\put(964,-3724){\makebox(0,0)[lb]{\smash{{\SetFigFont{14}{16.8}{\rmdefault}{\mddefault}{\updefault}{\color[rgb]{0,0,0}$w_{k+1}$}%
}}}}
\put(976,-4211){\makebox(0,0)[lb]{\smash{{\SetFigFont{14}{16.8}{\rmdefault}{\mddefault}{\updefault}{\color[rgb]{0,0,0}$w_{k+2}$}%
}}}}
\put(926,-6561){\makebox(0,0)[lb]{\smash{{\SetFigFont{14}{16.8}{\rmdefault}{\mddefault}{\updefault}{\color[rgb]{0,0,0}$w_{k+m}$}%
}}}}
\put(926,-3011){\makebox(0,0)[lb]{\smash{{\SetFigFont{14}{16.8}{\rmdefault}{\mddefault}{\updefault}{\color[rgb]{0,1,0}$w_k$}%
}}}}
\end{picture}%